\documentclass[journal, twocolumn]{IEEEtran}

\usepackage[pdftex]{graphicx}
\usepackage[cmex10]{amsmath}
\usepackage{amsfonts}
\usepackage{amssymb}
\usepackage{acronym}
\usepackage{color}
\usepackage{xcolor}
\usepackage[hidelinks]{hyperref}
\usepackage{comment}
\usepackage{subcaption}
\usepackage{threeparttable}
\usepackage{booktabs}
\usepackage{amsthm}
\usepackage{upgreek}
\usepackage{multirow}
\usepackage{cite}
\usepackage{scrextend}
\addtokomafont{labelinglabel}{\sffamily\bfseries}

\newcommand{\bfb}[1]{\boldsymbol{\rm #1}}
\newcommand{\bfg}[1]{\boldsymbol{#1}}
\newcommand{\nx}{\nu}
\newcommand{\ny}{\mu}

\newcommand{\xpc}{\xi}
\newcommand{\xs}{\ensuremath{\mbox{$\bfb {x}$}}}

\newcommand{\jac}[2]{\bfg{{#1}}_{\hspace{-0.2mm}#2}}

\newcommand{\jj}{\jmath}

\newcommand{\T}{^{\intercal}}
\newcommand{\AS}{\bfb A_{\rm s}}
\newcommand{\vS}{\bfb v}
\newcommand{\scsym}[1]{\scriptscriptstyle \rm #1}

\usepackage{array}
\newcolumntype{P}[1]{>{\centering\arraybackslash}p{#1}}

\theoremstyle{definition}
\newtheorem*{theorem*}{Theorem}

\acrodef{fem}[FEM]{Forward Euler Method}
\acrodef{hm}[HM]{Heun's Method}
\acrodef{tm}[TM]{Trapezoidal Method}
\acrodef{pc}[PC]{Predictor-Corrector}
\acrodef{psa}[PSA]{Partitioned-Solution Approach}
\acrodef{sssa}[SSSA]{Small-Signal Stability Analysis}

\acrodef{aiits}[AIITS]{All-Island Irish Transmission System}
\acrodef{bdf2}[BDF2]{2-step Backward Differentiation Formula}
\acrodef{bem}[BEM]{Backward Euler Method}
\acrodef{dae}[DAE]{Differential-Algebraic Equation}
\acrodef{dde}[DDE]{Delay Differential Equation}
\acrodef{sg}[SG]{Synchronous Generator}

\acrodef{ode}[ODE]{Ordinary Differential Equation}
\acrodef{tds}[TDS]{Time-Domain Simulation}
\acrodef{2sdirk}[2S-DIRK]{2-Stage Diagonally Implicit Runge-Kutta}
\acrodef{avr}[AVR]{Automatic Voltage Regulator}
\acrodef{tg}[TG]{Turbine Governor}
\acrodef{pss}[PSS]{Power System Stabilizer}

\definecolor{ethblau}{RGB}{33,92,175}

\title{Unified Numerical Stability and Accuracy Analysis of the Partitioned-Solution Approach} 

\author{%
  Georgios Tzounas, {\em IEEE Member},
  and Gabriela Hug, {\em IEEE Senior Member}
 \thanks{This work is supported by the Swiss National Science Foundation under NCCR Automation (grant no.~51NF40 18054).
  \textit{(Corresponding author:
Georgios Tzounas.)}}
  \thanks{The authors are with ETH Z{\"u}rich, 8092 Z{\"u}rich, Switzerland.  E-mails:
    \{georgios.tzounas, hug\}@eeh.ee.ethz.ch.}%
  }

\begin{document}

\maketitle \pagestyle{plain} \IEEEpeerreviewmaketitle

\begin{abstract} 

This paper focuses on the \ac{psa} employed for the \ac{tds} of dynamic power system models.
In \ac{psa}, differential equations are solved at each step of the \ac{tds} for state variables, whereas algebraic equations are solved separately. 
The goal of this paper is to propose a novel, matrix-pencil based technique to study numerical stability and accuracy of \ac{psa} in a unified way. 
The proposed technique quantifies the numerical deformation that \ac{psa}-based methods introduce to the dynamics of the power system model, and allows estimating useful
upper time step bounds that achieve prescribed simulation accuracy criteria.  The family of \ac{pc} methods, which is commonly applied in practical implementations of \ac{psa}, is utilized to illustrate the proposed technique.  Simulations are carried out on the IEEE 39-bus system, as well as on a 
1479-bus model of the \ac{aiits}.
\end{abstract}
%

\section{Introduction}
\label{sec:intro}

\subsection{Motivation}


Despite recent advances in non-linear stability theory and computational methods for Lyapunov functions, e.g.~see \cite{giesl2015review,abate2020formal,anderson2015advances}, the most successful method to assess the stability of electric power systems under large disturbances to date is running numerical \acfp{tds}.  Nevertheless, fast and accurate stability assessment through \acp{tds} is not a trivial problem, due to the complexity of dynamics and the large scale of power systems.  These challenges will be further exacerbated in the future 
for a variety of reasons, including the increasing penetration of converter-based resources, the increase in electric energy demand, etc.
The focus of this paper is on the numerical stability and accuracy of a \ac{tds} approach that is being widely utilized both in industry and academia as a means of providing fast calculations in dynamic power system studies, namely the partitioned-solution approach \cite{machowski2020power}.

In general, a \ac{tds} consists in employing a proper numerical scheme to compute the time-response of a given dynamic model under known initial conditions.  In power system short-term stability analysis, the  dynamic model of the system is conventionally formulated as a set of non-linear \acp{dae}, as follows \cite{kundur:94}:
\begin{equation}
  \begin{aligned}
    \label{eq:dae}
    {\bfg x}'(t)  &= \bfg f( \bfg  x(t) , f\bfg y (t) ) \, , \\
    \bfg 0 &=  \bfg  g( \bfg x(t), \bfg y (t) ) \, .
  \end{aligned}    
\end{equation}
In \eqref{eq:dae}, $\bfg x(t): [0,\infty) \rightarrow \mathbb{R}^{\nx}$ are the states 
of 
dynamic devices connected to the power network (like conventional and converter-based generation, dynamic loads, automatic controllers, etc.); and $\bfg y(t): [0,\infty) \rightarrow \mathbb{R}^{\ny}$ are algebraic variables (like bus voltage magnitudes and phase angles, bus power injections, auxiliary variables that define control set-points, etc.); $\bfg f : \mathbb{R}^{\nx+\ny} \rightarrow \mathbb{R}^{\nx}$ and $\bfg g : \mathbb{R}^{\nx+\ny} \rightarrow \mathbb{R}^{\ny}$ are column vectors of non-linear functions that define the differential and algebraic equations, respectively; $\bfg 0$ is the zero matrix of proper dimensions.  Finally, discrete dynamics (like actions of control limiters, switching devices, etc.) do not appear in \eqref{eq:dae} because they are modeled implicitly through \textit{if-then rules}, i.e.,~each discontinuous change leads to a jump from \eqref{eq:dae} to another \ac{dae} set of the same form \cite{hiskens2004power}.

Various numerical schemes have been reported in the literature for the solution of \eqref{eq:dae}.  An important component of each scheme is the strategy it adopts to handle the equations, in particular whether differential and algebraic equations are solved together at once, or their solution is in some way alternated.  The former approach is referred to as \textit{simultaneous}, while the latter as \textit{partitioned}
\cite{stott:1979,machowski2020power}.  A second important component is the type of 
integration method employed, i.e.~\textit{implicit} or \textit{explicit}.  With this regard, the main advantages of implicit over explicit methods are their properties of numerical stability and their ability to deal with system \textit{stiffness} \cite{kundur:94}. Yet, these advantages are moderated by the need for frequent full Jacobian matrix factorizations, which increase the computational cost of \ac{tds} per time step.  On the other hand, explicit methods are computationally cheaper per time step; however, they give rise to large errors and numerical instability unless a small-enough time step is used which, in turn, also increases the total computational time.  Finally, a third important component is whether the solution is obtained by exploiting parallelization algorithms, e.g.~see \cite{aristidou2013schur,gurrala2015parareal,jin2017comparative}.

The problem of choosing a good numerical integration scheme is one of finding a good compromise between speed and accuracy/stability.  
In practice, implicit numerical methods are most commonly combined with the simultaneous 
approach.  Methods employed in this context include the \ac{tm}, backward Euler, Theta, diagonally implicit Runge-Kutta, Hammer-Hollingsworth, etc., see~\cite{2009noda, hh4:2013, stiffdecay,1995paserba}.  On the contrary, explicit methods are implemented with the \acf{psa}.  In this context, very simple methods such as the \ac{fem} are avoided due to their poor performance, with most tools choosing to combine an explicit method with some accuracy refinement technique, such as \acf{pc} iterations
\cite{kundur:94, machowski2020power}.  
However, even with accuracy refinement, \ac{psa}-based solutions are 
still prone to numerical issues,  which, once again, limits the ability to use large time steps. This has often driven efforts for the definition of device models that are numerically robust when combined with a given commercial \ac{psa}-based solver.  In this vein, we refer to
\cite{ramasubramanian2020positive,ramasubramanian2022parameterization} on the 
development and parameterization 
of 
power converter models for systems with low short-circuit strength.

The main goal of this paper is to propose a novel technique to study in a unified way the stability and precision of numerical integration of power system models with \ac{psa}.  
We note that classical stability characterization based on the response of scalar test equations is not a suitable approach to this aim, as it omits the differential and algebraic equations of the examined system.
In this regard, a technique to study \ac{psa} by taking into account the system's equations 
 is proposed in \cite{milanopsa}. Therein, the effect of \ac{psa} is seen as that of an 
 \textit{one-step delay} introduced to all algebraic variables that appear in the system's differential equations.  The main limitation of the technique in \cite{milanopsa} is that it neglects the effect of the integration method applied, which in turn prevents a precise estimation of the available margin until 
 numerical instability.
Finally, a matrix-based framework to study \ac{tds} factoring in the effect of both system dynamics and integration method was recently introduced by the first author of this paper in \cite{tzounas2022tdistab,tzounas2022tdistabd} without, however, studying the \ac{psa}.  In fact, the numerical methods and accuracy refinement techniques commonly employed in the implementation of \ac{psa}
can not be studied by means of the derivations in these works. 



\subsection{Contributions}


Given the limitations of existing techniques identified in the previous section, this paper 
proposes a novel, matrix pencil-based approach to study in a unified way the numerical stability and accuracy of \ac{psa} employed for the solution of the \acp{dae} that describe the evolution of power system dynamics.  

By employing the family of schemes most commonly used in practical implementations of \ac{psa}, i.e.~\ac{pc} methods, the paper quantifies the deformation introduced to the computed system dynamics, first, due to the application of the integration method \textit{per se} and, second, due to the mismatch between state and algebraic variables caused by the so called \textit{interfacing problem}. 
This information can be then used to find the numerical stability margin as well as to estimate useful upper time step bounds that achieve prescribed simulation accuracy criteria.
For \acf{hm}, the matrix pencils that need to be studied 
to this aim are fully derived, both for sparse and dense matrix calculations. 
Finally, the paper provides a comparison of the proposed technique with the one-step delay method described in \cite{milanopsa}.


\subsection{Organization}

The remainder of the paper is organized as follows.  Section~\ref{sec:preliminaries} provides preliminaries on \ac{tds} of power systems through implicit and explicit integration methods.  The \ac{psa} and its implementation using \ac{pc} methods are described in Section~\ref{sec:psa}.  Section~\ref{sec:stability}
provides the formulation of the proposed matrix pencil-based technique to assess the numerical stability and accuracy of \ac{psa}.  
The case studies are discussed in Section~\ref{sec:case} based on the IEEE 39-bus system, as well as on a realistic 1479-bus model of the \acf{aiits}.
Finally, conclusions are drawn and future work directions are outlined in Section~\ref{sec:conclusion}.

\section{Preliminaries}
\label{sec:preliminaries}

\subsection{Notation}

The notation adopted throughout the paper is as follows: $a$, $\bfg a$, $\bfb A$, denote, respectively, scalar, vector, and matrix, and $\bfb A\T$ denotes the matrix transpose; $\mathbb{N}$, $\mathbb{R}$, $\mathbb{C}$ are the sets of natural, real, and complex numbers, respectively, and $\jj$ is the unit imaginary number; 
$a^{(i)}$ denotes the $i$-th iteration;
$a(t)$ denotes a continuous-time quantity and
$a'(t)$ its first-order derivative;
$\mathcal{L}\{\cdot\}$ is the Laplace transform,
$s$ is the complex frequency in $S$-domain,
and $a(s)$ is a $S$-domain quantity; $a_{n+1}$ denotes a discrete-time quantity; $z$ is the complex frequency in $Z$-domain; $\Delta$ and $\nabla$ are the forward and backward difference operators.

\subsection{Numerical Integration Methods}

A numerical integration method applied to a \ac{dae} power system model can be formulated as a set of non-linear difference equations, the solution of which is an approximation of the time-domain response of \eqref{eq:dae}.  

\subsubsection{Implicit Methods}

Let's define:
\begin{equation}
\begin{aligned}
\label{eq:discrxy}
\bfg x_{n+1-\ell} &= \bfg x(t+(1-\ell) h) \, , \\
\bfg y_{n+1-\ell} &= \bfg y(t+(1-\ell) h) \, ,
\end{aligned}    
\end{equation}
where $\bfg x_{n+1-\ell} : \mathbb{N} \rightarrow \mathbb{R}^{\nx}$ and $\bfg y_{n+1-\ell} : \mathbb{N} \rightarrow \mathbb{R}^{\ny}$ denote the 
discretized state and algebraic variables, respectively, at time $t+(1-\ell) h$, with $l \in \mathbb{N}$.
Then, an implicit integration method for the solution of \eqref{eq:dae} can be described as follows:
\begin{equation}
\begin{aligned}
  \label{eq:tdi:implicit}
  \bfg 0 &= \bfg \phi(
  \bfg x_{n+1}, \bfg y_{n+1},
  \bfg x_{n}, \bfg y_{n}, 
  \ldots,
  \bfg x_{n-a}, \bfg y_{n-a}, h) \, , \\ 
  \bfg 0 &= \bfg \rho(
  \bfg x_{n+1}, \bfg y_{n+1},
  \bfg x_{n}, \bfg y_{n}, \ldots,
  \bfg x_{n-a}, \bfg y_{n-a}, h) \, , \\ 
\end{aligned}
\end{equation}
where $h$ is the simulation time step size;
$\bfg x_{n-a}, \bfg y_{n-a}$, with $a\in\mathbb{N}^*$, are included 
so that \eqref{eq:tdi:implicit} 
covers 
both Runge-Kutta and multistep methods.
Moreover, $\bfg\phi$, $\bfg\rho$ are column vectors of 
functions that depend on the differential and algebraic equations, respectively, as well as on the specific numerical method applied.  
For example, for the \ac{tm} we have:
\begin{equation}
\begin{aligned}
\label{eq:tm}
  \bfg 0 &= \bfg x_{n+1} - \bfg x_{n} - 0.5h
 [\bfg f(\bfg x_{n},\bfg y_{n}) + \bfg f(\bfg x_{n+1},\bfg y_{n+1})] \, , \\
 \bfg 0 &= h \bfg g(\bfg x_{n+1},\bfg y_{n+1}) \, .
\end{aligned}
\end{equation}

Solving 
\eqref{eq:tdi:implicit} at every time step consists in calculating the solution $(\bfg x_{n+1}, \bfg y_{n+1}) := [\bfg x_{n+1}\T, \bfg y_{n+1}\T]\T$ (where $\T$ indicates the transpose), which in this case is typically obtained through the simultaneous-solution approach, i.e.~by updating together the state and algebraic variables of the system.  This is done through Newton iterations, i.e.~by solving the problem:
\begin{equation}
\label{eq:newton}
\bfb J
\begin{bmatrix}
\Delta \bfg x_{n+1}^{(i)} \\
\Delta \bfg y_{n+1}^{(i)} \\
\end{bmatrix}
= 
-
\begin{bmatrix}
\bfg \phi^{(i)} \\
\bfg \rho^{(i)} \\
\end{bmatrix} 
, \ \ \text{where} \ \
\bfb J =
\begin{bmatrix}
\jac{\phi}{x}^{(i)} & \jac{\phi}{y}^{(i)} \\
\jac{\rho}{x}^{(i)} & \jac{\rho}{y}^{(i)}
\end{bmatrix} ,
\end{equation}
with $\jac{\phi}{x}^{(i)}={\partial {\bfg \phi^{(i)}}}/{\partial \bfg x^{(i)}_{n+1}}$ (similarly for the rest); and then determining the new values as:
%
\begin{equation}
\begin{aligned}
\label{eq:update}
\bfg x_{n+1}^{(i+1)} &=
\bfg x_{n+1}^{(i)} + \Delta \bfg x_{n+1}^{(i)} \, , \\
\bfg y_{n+1}^{(i+1)} &=
\bfg y_{n+1}^{(i)} + \Delta \bfg y_{n+1}^{(i)} \, .
\end{aligned}
\end{equation}
The procedure above is repeated until convergence, i.e.~until $||\bfg x_{n+1}^{(i+1)} - \bfg x_{n+1}^{(i)}|| < \epsilon$ and $||\bfg y_{n+1}^{(i+1)} - \bfg y_{n+1}^{(i)}|| < \epsilon$, where $\epsilon>0$ is a given tolerance.  
The most costly step in this process is the factorization of the Jacobian matrix $\bfb J$, which is required for the solution of \eqref{eq:newton}.  
Since $\bfb J$ is a general, non-symmetric, sparse matrix, its factorization is most commonly achieved through symbolic and numeric LU decomposition \cite{davis2006direct}.
Moreover, to speed up the solution, a \textit{(very) dishonest} Newton method is often used, where $\bfb J$ is factorized only once per (multiple) time step(s).  We note that a dishonest scheme is, ultimately, a compromise, as it can speed up calculations but also sacrifices convergence and may lead to numerical issues, such as infinite-cycling.



\subsubsection{Explicit Methods}

In contrast to implicit integration methods, explicit methods calculate the new state vector $\bfg x_{n+1}$ at each step without the need to factorize $\bfb J$.  In general, an explicit numerical method for the solution of \eqref{eq:dae} can be described as follows:
\begin{equation}
\begin{aligned}
\label{eq:tdi:explicit}
\bfg x_{n+1} &= 
\bfg \psi(\bfg x_{n}, \bfg y_{n}, \ldots,
\bfg x_{n-a}, \bfg y_{n-a}, h) \, , \\ 
\bfg 0 &= h\bfg g(\bfg x_{n+1},\bfg y_{n+1}) \, .
\end{aligned}
\end{equation}
For example, if \ac{fem} is employed, we have:
\begin{align}
\label{eq:fem:f}
\bfg x_{n+1} &= 
\bfg x_{n} + {h} 
\bfg f(\bfg x_{n},\bfg y_{n}) \, , \\
\label{eq:fem:g}
\bfg 0 &= h\bfg g(\bfg x_{n+1},\bfg y_{n+1}) \, .
\end{align}
Solution of \eqref{eq:tdi:explicit} is usually obtained through \ac{psa}.  Since \ac{psa} is the main focus of this paper, we discuss it in detail in 
Section~\ref{sec:psa}.

\subsection{Classical Stability Characterization}

In this section, we briefly recall the classical approach to study the stability of a numerical integration method. 
The stability region of a given numerical integration method 
is conventionally characterized by evaluating its response when applied to a scalar, linear test differential equation:
\begin{equation}\label{eq:test}
  x'(t)=\lambda x(t) \, , 
  \ \ \lambda \in \mathbb{C} \, .
\end{equation}
%
Let's apply an integration method to \eqref{eq:test} so that:
\begin{equation}\label{eq:stabfunc}
  x_{n+1}=\mathcal{R}(\lambda h) x_{n} \, .
\end{equation}
Then, $\mathcal{R}(\lambda h)$ is called the method's \textit{growth function} and 
the stability region of the method is defined by the set:
%
\begin{equation}
  \label{eq:test:region}
  \{ \lambda \in \mathbb{C} \, :
  \  |\mathcal{R}(\lambda h)|<1 \}
  \, .
\end{equation}
Considering, for example, the \ac{tm}, we have:
%
\begin{equation}\label{eq:itm:test}
  x_{n+1}=x_{n} + 0.5h\lambda x_{n} + 0.5h\lambda x_{n}
  \, ,
\end{equation}
which can be equivalently written in the form of \eqref{eq:stabfunc}, where:
\begin{equation}\label{eq:itm:stability}
  \mathcal{R}(\lambda h) =
  ({1+0.5\lambda h})/({1-0.5\lambda h}) \, . 
\end{equation}
From \eqref{eq:test:region}, \eqref{eq:itm:stability}, we deduce that the region of stability of the \ac{tm} is the left half of the $S$-plane.

The main limitation of the above approach is that it neglects the 
equations of the examined 
model.  Hence, it can be only used for certain methods and only qualitatively, and thus it is inherently unsuitable for accuracy analysis.

\section{The Partitioned-Solution Approach}
\label{sec:psa}

In this section, we provide a description of \ac{psa} utilized for the numerical integration of a \ac{dae} power system model in the form of \eqref{eq:dae}.

\subsection{Forward Euler Method}
\label{sec:fem}

Let us initially consider for the sake of illustration that \ac{fem} is employed for the solution of \eqref{eq:dae}, and that $(\bfg x_n, \bfg y_n)$ is known for some $n$.  Then, the idea of 
\ac{psa} is to obtain the solution $(\bfg x_{n+1}, \bfg y_{n+1})$ in an alternating fashion.  That is, first, $(\bfg x_n, \bfg y_n)$ is used in \eqref{eq:fem:f} to evaluate $\bfg f$ and then calculate $\bfg x_{n+1}$.  Second, using the calculated $\bfg x_{n+1}$, \eqref{eq:fem:g} is solved for $\bfg y_{n+1}$, e.g.~through Newton iterations, which in this case read as follows: 
\begin{equation}
\begin{aligned}
\label{eq:fem:newton}
\jac{g}{y}^{(i)}
\Delta \bfg y_{n+1}^{(i)} &= - h\bfg g^{(i)} \, , \\
\nonumber
\bfg y_{n+1}^{(i+1)} &=
\bfg y_{n+1}^{(i)} + \Delta \bfg y_{n+1}^{(i)} \, .
\end{aligned}
\end{equation}
\ac{fem} is the simplest and most intuitive among all integration methods.  However, it also shows a poor performance and hence, \ac{tds} routines of 
software tools that adopt \ac{psa} prefer to rely on more robust and complex schemes.  With no loss of generality, in the following we further discuss \ac{psa} by considering \acfp{pc} schemes, i.e.~the family of methods most commonly used in practical implementations of \ac{psa} \cite{stott:1979, milanopsa}.

\subsection{Predictor-Corrector (PC) Iterations}
\label{sec:predcorr}

\ac{pc} schemes combine the application of explicit and implicit numerical methods.  The 
steps of a generic \ac{pc} scheme employed for the solution of \eqref{eq:dae} can be described as follows:

\vspace{0.6mm}
\textit{Predictor}:  The predictor implements an explicit method that provides an initial estimation ($\bfg\xpc^{(0)}_{n+1}$) of $\bfg x_{n+1}$, as follows:
\begin{align}\label{eq:pred}
  \bfg\xpc^{(0)}_{n+1} &= 
  \bfg \psi(\bfg x_{n}, \bfg y_{n}, \ldots,
  \bfg x_{n-a}, \bfg y_{n-a}, h) \, .
\end{align}

\vspace{0.6mm}
\textit{Corrector}: The corrector employs an implicit method to refine the accuracy of the estimation of $\bfg x_{n+1}$.  The $i$-th corrector step has the form:
\begin{equation}\label{eq:corr}
\hspace{-1.8mm}
  \bfg 0  = 
  \bfg \phi( \bfg\xpc^{(i)}_{n+1}, 
  \bfg\xpc^{(i-1)}_{n+1},
  \bfg y_{n+1},
  \bfg x_{n}, \bfg y_{n}, \ldots,
  \bfg x_{n-a}, \bfg y_{n-a}, h) \, ,
\end{equation}
where $i \in \mathbb{N}^* : i\leq r$, which means that the accuracy is improved by iteration (typical values are $r=1$ or $r=2$).  

Then, $\bfg x_{n+1}$ is obtained simply as:
\begin{equation}\label{eq:corr:xt}
\bfg x_{n+1} = \bfg \xpc^{(r)}_{n+1} \, .
\end{equation}
The scheme is completed considering the algebraic equations:
\begin{equation}
\label{eq:pc:g}
 \bfg 0 = h \bfg g(\bfg x_{n+1},\bfg y_{n+1}) \, .
\end{equation}
To solve for $(\bfg x_{n+1},\bfg y_{n+1})$
for a generic scheme in the form of \eqref{eq:pred}-\eqref{eq:pc:g},\footnote{Note that there exist multiple ways in which the solution of state/algebraic variables can be updated and thus formulation \eqref{eq:pred}-\eqref{eq:pc:g} is not unique. E.g.,~an alternative is 
to solve algebraic variables after each corrector step
(instead of only in the end).  Comparing all possible formulations is not in the scope of this paper, yet, solving more often for algebraic variables is expected to improve accuracy but also increase the computational cost. 
}
the basic procedure of \ac{psa} is in analogy
to the one discussed in the beginning of Section~\ref{sec:fem} for \ac{fem}.  That is, $\bfg x_{n+1}$ is first calculated from \eqref{eq:pred}-\eqref{eq:corr:xt}. Then, $\bfg y_{n+1}$ is obtained from \eqref{eq:pc:g}.  Nevertheless, notice that, in contrast to \ac{fem}, using \eqref{eq:pred}-\eqref{eq:corr:xt} to calculate $\bfg x_{n+1}$ requires knowing $\bfg y_{n+1}$ (see \eqref{eq:corr}), which is though not available yet.  The problem of handling the unknown $\bfg y_{n+1}$ in \eqref{eq:corr} is an important problem of \ac{psa}, commonly referred to as the \textit{interfacing problem}.

\subsection{Interface Errors}
\label{sec:if}

The interfacing problem 
is a consequence of the inherent coupling of the \acp{dae} used to describe the dynamic behavior of power systems.  In general, the way that software tools deal with this problem is by adopting an approximation, where $\bfg y_{n+1}$ is substituted by a known value, say $\bfg y_{\rm int}$.  For \ac{pc} schemes, this leads to a change of \eqref{eq:corr} into:
\begin{equation}
  \label{eq:corr:if}
  \bfg 0 = \bfg \phi( \bfg\xpc^{(i)}_{n+1}, 
  \bfg\xpc^{(i-1)}_{n+1}, \bfg y_{\rm int},
  \bfg x_{n}, \bfg y_{n}, \ldots,
  \bfg x_{n-a}, \bfg y_{n-a}, h) \, .
\end{equation}

Equation \eqref{eq:corr:if} describes the corrector in a general way without imposing a specific strategy to deal with the interfacing problem.  Apparently, once the interfacing strategy is fixed, $\bfg y_{\rm int}$ is fully defined.  In this regard, the simplest and computationally most efficient interfacing strategy is extrapolation, i.e.~the value at the previous step ($\bfg y_n$) is used instead of $\bfg y_{n+1}$.  However, using $\bfg y_{\rm int}=\bfg y_{n}$ in \eqref{eq:corr:if} introduces a mismatch between state and algebraic variables and thus an error to the numerical solution of the system, known as the \textit{interface error}.  One approach to reduce the interface error is comparing $\bfg y_{\rm int}$ with $\bfg y_{n+1}$ once the latter has been computed.  If the difference is higher than a given threshold, then the integration of the step is repeated until an acceptable tolerance is achieved, or equivalently, until $\bfg y_{\rm int} \approx \bfg y_{n+1}$.

\section{Pencil-Based Stability Analysis}
\label{sec:stability}

In this section, we first define the stiffness and \ac{sssa} of \eqref{eq:dae}, 
and then proceed to describe the proposed \ac{sssa}-based technique to study \ac{psa}. 


\subsection{Model Stiffness and SSSA}
\label{sec:dae:sssa}

Power system models are known to be stiff, i.e.~the time constants that define the differential equations of the model span multiple time scales.  The stiffness of the \ac{dae} system \eqref{eq:dae} is typically measured by the ratio between the largest and smallest eigenvalues of the corresponding small-signal model.
Assume that an equilibrium $(\bfg x_o, \bfg y_o)$ 
of \eqref{eq:dae} is known.  At the equilibrium, we have $\bfg 0 = \bfg f(\bfg x_o , \bfg y_o )$, 
$\bfg 0 = \bfg g(\bfg x_o , \bfg y_o )$.  For sufficient small disturbances and for the purpose of applying well-known results from linear stability theory, \eqref{eq:dae} can be linearized in a neighborhood of $(\bfg x_o, \bfg y_o)$ as follows:
\begin{equation}
  \begin{aligned}
    \label{eq:dae:lin}
    \tilde {\bfg x}'(t)  &= 
    \jac{f}{x} \tilde {\bfg x}(t) + \jac{f}{y} \tilde {\bfg y}(t) \, , \\
    \bfg 0 &=\jac{g}{x} \tilde {\bfg x}(t) + \jac{g}{y} \tilde {\bfg y}(t) \, , 
  \end{aligned}    
\end{equation}
where $\tilde {\bfg x}(t) = \bfg x(t) - \bfg x_o$, $\tilde {\bfg y}(t) = \bfg y(t) - \bfg y_o$; and $\jac{f}{x}$, $\jac{f}{y}$, $\jac{g}{x}$, $\jac{g}{y}$ are Jacobian matrices evaluated at $(\bfg x_o, \bfg y_o)$.  Applying the Laplace transform to \eqref{eq:dae:lin} and omitting for simplicity the time-dependency in the notation, we have:
\begin{equation}\label{eq:dae:lin:lap}
\begin{bmatrix}
s \mathcal{L} \{\tilde {\bfg x} \} - \tilde {\bfg x}(0) \\
{\bfg 0} \\
\end{bmatrix}
= 
\begin{bmatrix}
\jac{f}{x} & \jac{f}{y} \\
\jac{g}{x} & \jac{g}{y} \\
\end{bmatrix}
\begin{bmatrix}
 \mathcal{L} \{\tilde {\bfg x} \}  \\
 \mathcal{L} \{\tilde {\bfg y} \} \\
\end{bmatrix}
 ,
\end{equation}
where $s$ is a complex frequency in the $S$-plane. Equivalently:
\begin{equation}
\label{eq:dae:lin:lap2}
(s \bfb E - \bfb A) 
\mathcal{L} 
\left \{
\begin{bmatrix}
\tilde {\bfg x} \\
\tilde {\bfg y} \\
\end{bmatrix}
\right \}
= \bfb E 
\begin{bmatrix}
\tilde {\bfg x}(0) \\
 {\bfg 0} \\
\end{bmatrix} ,
\end{equation}
with
\begin{equation}\label{eq:dae:EA}
\bfb E = 
  \begin{bmatrix}
    \bfb I & \bfg 0 \\
    \bfg 0 & \bfg 0 \\
  \end{bmatrix} ,  \ \ 
\bfb A = 
\begin{bmatrix}
\jac{f}{x} & \jac{f}{y} \\
\jac{g}{x} & \jac{g}{y} \\
\end{bmatrix}
,
\end{equation}
where $\bfb I$ denotes the identity matrix of proper dimensions.  

Then, the polynomial matrix
\begin{equation}\label{eq:dae:pencil:sparse}
\bfb P(s)=s\bfb E-\bfb A \, , 
\end{equation}
is called the \textit{matrix pencil} of \eqref{eq:dae:lin} and plays an important role in the study of the system.  
The eigenvalues of this pencil
can be obtained from the solution of the algebraic problem \cite{book:eigenvalue}:
\begin{equation}\label{eq:ep}
\bfb P(s)  \bfb v = \bfg 0  \, ,
\end{equation}
with $\bfb v \in \mathbb{C}^{(\nx+\ny)\times 1}$.
Then, system \eqref{eq:dae:lin} is 
stable if for every eigenvalue $s_i$ of $\bfb P(s)$, $\Re(s_i)<0$.
%
Let the system be stable and $s^{\max}$, $s^{\min}$ be the eigenvalues of \eqref{eq:dae:pencil:sparse} with the largest and smallest magnitudes, i.e.~$s^{\max}=\max{|s_i|}$,
$s^{\min}=\min{|s_i|}$, $\forall s_i$, then the stiffness ratio of \eqref{eq:dae} can be defined as follows:
\begin{equation}\label{eq:stiffness}
\mathcal{S} = |s^{\max}|/|s^{\min}| \, .
\end{equation}

Alternatively to \eqref{eq:dae:pencil:sparse}, algebraic variables can be eliminated from 
\eqref{eq:dae:lin}
provided that $\bfg g_y$ is 
non-singular, 
in which case the following pencil can be equivalently used in \eqref{eq:ep}:
\begin{equation}\label{eq:dae:pencil:dense}
\bfb P(s) = s\bfb I-\AS  \, , 
\end{equation}
where $\AS =\jac{f}{x}-\jac{f}{y}\jac{g}{y}^{-1}\jac{g}{x}$; and $\vS \in \mathbb{C}^{\nx\times 1}$.

We have two notes on \eqref{eq:dae:pencil:sparse}, \eqref{eq:dae:pencil:dense}.  

First, the two pencils have the same finite eigenvalues, yet  \eqref{eq:dae:pencil:sparse} has also the infinite eigenvalue with multiplicity $\ny$.  
Recall that finite eigenvalues are the zeros of the polynomial $\det(s\bfb E - \bfb A) = p(s)$, which is of order $\nx$.  On the other hand, the existence of the infinite eigenvalue is seen by means of the 
problem $(s\bfb E - \bfb A) \bfb v = \bfg 0 $ rewritten in the reciprocal form
$(\bfb E - s^{-1}\bfb A) \bfb u = \bfg 0 $, where 
$\bfb u \in \mathbb{C}^{\nx\times 1}$.
Since $\bfb E$ is singular, there exists a null vector $\bfb u$ such that $\bfb E \bfb u = \bfb 0$, and thus
$s^{-1}\bfb A \bfb u = \bfg 0 $,
so that $\bfb u$ is an eigenvector of the reciprocal problem corresponding to the eigenvalue $s^{-1}=0$, or $s \rightarrow \infty$~\cite{moebius}.

Second, $\bfb E$ and $\bfb A$ are sparse matrices, whereas $\AS$ is dense.  This makes a difference in the numerical calculations.  For small/medium size systems, \eqref{eq:dae:pencil:dense} is preferred and the solution is obtained with the QR algorithm \cite{magma,app10217592}.  On the other hand, \eqref{eq:dae:pencil:sparse} is preferred for large-scale systems and when it suffices to find only a subset of the spectrum.  The solution in this case is found with a sparse algorithm that can handle asymmetric matrices, such as Krylov-Schur or contour integral-based methods \cite{app10217592,li2015ks,li2016cirr}.

\subsection{Proposed Numerical Stability and Accuracy Assessment}\label{sec:adams:sssa}

In this section, we describe the proposed approach.
To this aim, we consider a common \ac{pc} method, namely Heun's method (or modified Euler method), variants of which are used in commercial tools, including \cite{psse_pag2}, \cite{ge:pslf}.
\acf{hm} is a combination of \ac{fem} and \ac{tm}.  When applied to \eqref{eq:dae}, the method reads as follows:
\begin{equation}
\begin{aligned}
\label{eq:mem}
\bfg \xpc^{(0)}_{n+1} &= \bfg x_{n} + h 
\bfg f(\bfg x_{n},\bfg y_{n}) \, , \\
\bfg \xpc^{(i)}_{n+1} &= \bfg x_{n}
+ 0.5h
\bfg f(\bfg x_{n},\bfg y_{n})
+  0.5h \bfg f (\bfg \xpc^{(i-1)}_{n+1}, 
\bfg y_{\rm int}) \, , \\
\bfg x_{n+1} &= \bfg \xpc^{(r)}_{n+1} \, , \\
\bfg 0 &= h \bfg g(\bfg x_{n+1},\bfg y_{n+1}) \, .
\end{aligned}
\end{equation}
Notice that, due to $\bfg\xpc^{(i-1)}_{n+1}$ being always available, the corrector can be expressed in an explicit form, without the need to resort to further calculations.


Let $(\bfg x_o,\bfg y_o)$ be an equilibrium of \eqref{eq:dae}.  Then, $(\bfg x_o,\bfg y_o)$ is also a fixed point of \eqref{eq:mem}.
Linearizing \eqref{eq:mem} at this point gives:
\begin{align}
\label{eq:mem:pred:lin}
\tilde {\bfg \xpc}^{(0)}_{n+1} &= \tilde {\bfg x}_{n} + h (\jac{f}{x} \tilde {\bfg x}_{n} + \jac{f}{y} \tilde {\bfg y}_{n} ) \, , \\
\label{eq:mem:cori:lin}
\tilde {\bfg \xpc}^{(i)}_{n+1} &= \tilde {\bfg x}_{n} + \frac{h}{2} [\jac{f}{x}\tilde {\bfg x}_{n} + \jac{f}{x} \tilde {\bfg \xpc}^{(i-1)}_{n+1} + \jac{f}{y} (\tilde {\bfg y}_{n} + \tilde {\bfg y}_{\rm int} ) ] \, , \\
\label{eq:corr:xt:lin}
\tilde {\bfg x}_{n+1} &= \tilde {\bfg \xpc}^{(r)}_{n+1} \, , \\
\label{eq:pc:g:lin}
\bfg 0 &= \jac{g}{x} \tilde {\bfg x}_{n+1}+
 \jac{g}{y} \tilde {\bfg y}_{n+1} \, .
\end{align}
Equivalently, \eqref{eq:mem:pred:lin}-\eqref{eq:pc:g:lin} can be rewritten as:
\begin{align}\nonumber
\begin{bmatrix}
\bfb I & \hspace{-1mm} \bfg 0\\
\jac{g}{x}& \hspace{-1mm} \jac{g}{y} \\
\end{bmatrix}
\hspace{-1mm}
\begin{bmatrix}
\tilde {\bfg x}_{n+1} \\
\tilde {\bfg y}_{n+1} \\
\end{bmatrix}
=&
\begin{bmatrix}
\bfb I + h \bfb C_r \jac{f}{x} & 
\hspace{-1mm}
( h \bfb C_{r} - \frac{h}{2} \bfb C_{r-1} )\jac{f}{y} \\
\bfg 0 & \bfg 0\\
\end{bmatrix}
\begin{bmatrix}
\tilde {\bfg x}_{n} \\
\tilde {\bfg y}_{n} \\
\end{bmatrix}
\\ \label{eq:mem:predcor:lin}
&+
\begin{bmatrix}
\frac{h}{2} \bfb C_{r-1}
\jac{f}{y} \tilde {\bfg y}_{\rm int} \\
\bfg 0 \\
\end{bmatrix} \, ,
\end{align}
\begin{equation}\label{eq:mem:cr}
\hspace{-2cm} 
\text{where}
\hspace{1.5cm}
\bfb C_r = \sum_{j=0}^r \left(\frac{h}{2} \jac{f}{x}\right )^j , \ \ r \in \mathbb{N}^* \, .
\end{equation}
The proof of \eqref{eq:mem:predcor:lin} is provided in Section~\ref{sec:proof:mem} of the Appendix.  

Expression \eqref{eq:mem:predcor:lin} can be further simplified once a specific interfacing strategy is adopted.  As discussed in Section~\ref{sec:if}, the most common interfacing strategy is extrapolation.  
In this case, ${\bfg y}_{\rm int} = \bfg y_{n}$ and \eqref{eq:mem:predcor:lin} corresponds to a linear system of difference equations with matrix pencil: 
\begin{align}
\label{eq:mem:pencil:ext}
{\bfb P}_{\scsym PC}(z) = z 
\begin{bmatrix}
\bfb I & \bfg 0\\
\jac{g}{x}& \jac{g}{y} \\
\end{bmatrix} 
-  \begin{bmatrix}
\bfb I + h\bfb C_r\jac{f}{x} & h\bfb C_r\jac{f}{y}  \\
\bfg 0  &  \bfb 0 \\
\end{bmatrix}
 ,
\end{align}
where $z$ is a complex frequency in the $Z$-plane.  
    
The eigenvalues of the pencil ${\bfb P}_{\scsym PC}(z)$ provide insights into the small-signal dynamics of the \ac{dae} system \eqref{eq:dae} approximated by \ac{hm} when interfacing is achieved by extrapolation.  Recall from Section~\ref{sec:dae:sssa} that the actual small-signal dynamics of \eqref{eq:dae} are represented by the eigenvalues of ${\bfb P}(s)$.  Thus, by comparing the eigenvalues of the two pencils we can determine \textit{how much} the \ac{hm} 
numerically deforms the dynamic modes of the power system model.  Moreover, by increasing the time step we can determine the corresponding numerical stability margin.  
In particular, provided that the small-signal power system model is stable, 
numerical stability is maintained for all time step sizes which lead to all eigenvalues of ${\bfb P}_{\scsym PC}(z)$ having magnitudes smaller than 1.
Note that in order for the dynamics of ${\bfb P}(s)$ and ${\bfb P}_{\scsym PC}(z)$ to be comparable, they need to be referred to the same domain.  In this paper, eigenvalues are always referred to the $S$-plane.  With this regard, recall that a complex frequency $z$ in the $Z$-plane is mapped to the $S$-plane through the relationship:
\begin{equation}\label{eq:z2s}
\hat 
s = {{\rm log}(z)}/{h} \, .
\end{equation}

Let us now assume that the interface error in some manner is eliminated or is negligible. 
In this case, $\tilde {\bfg y}_{\rm int} = \bfg y_{n+1}$ and \eqref{eq:mem:predcor:lin} corresponds to a linear system of difference equations with matrix pencil:
\begin{equation}
\hspace{-0.15cm}
\label{eq:mem:pencil:rig}
{\bfb {\bar P}}_{\scsym PC}(z) = 
z
\hspace{-0.7mm}
\begin{bmatrix}
\bfb I & 
\hspace{-2.2mm} 
-\bfb B \jac{f}{y}\\
\jac{g}{x} & \jac{g}{y} \\
\end{bmatrix}  
\hspace{-0.8mm}
- 
\hspace{-0.7mm}
\begin{bmatrix}
\bfb I 
\hspace{-0.2mm}
+ 
\hspace{-0.2mm}
h \bfb C_r \jac{f}{x} & 
\hspace{-1.2mm}
( h \bfb C_{r} 
\hspace{-0.7mm}
- \bfb B)\jac{f}{y} \\
\bfb 0 & \bfb 0 \\
\end{bmatrix}  
\hspace{-0.2mm} ,
\end{equation}
where $\bfb B = \frac{h}{2} \bfb C_{r-1}$.  Then, the eigenvalues of the pencil ${\bfb {\bar P}}_{\scsym PC}(z)$ represent the small-signal dynamics of \eqref{eq:dae} approximated by \ac{hm} and when perfect interfacing is assumed.
As a result, by comparing the dynamics of the ${\bfb {P}}_{\scsym PC}(z)$ with those of ${\bfb {\bar P}}_{\scsym PC}(z)$, we can quantify the amount of deformation of system dynamics that is explicitly due to the interface error.
This is further discussed in the case study presented in Section~\ref{sec:case}.


The following remarks 
are 
relevant:

\begin{itemize}
    \item The proposed approach can be used to 
estimate the time step size that allows achieving prescribed simulation accuracy criteria, taking into account the deformation of the system's dynamics.  
This is in contrast to current practice of PSA-based commercial software packages, which rely on purely empirical rules to estimate the maximum admissible time step, based on
long experience with the simulation of conventional, synchronous machine-based power systems, see \cite{psse_pag2}.
For example, it is common to fix the time step to 1 order of magnitude lower than the 
smallest time constant in the process being simulated, e.g.~see \cite{psse_pag2}. 
Yet, the efficacy of such heuristics is currently challenged by the gradual shift to converter-based power systems.  In particular, it has been seen that state-of-art converter models often lead to numerical instability when connected to a power network with low short circuit strength \cite{ramasubramanian2020positive,ramasubramanian2022parameterization}. 
On the contrary, the proposed technique 
can be applied to every system and is thus
device-independent, 
and this will be also true for future power system technologies. 


\item Different \ac{psa}-based methods approximate the 
solution of the \acp{dae} with different accuracies when applied under the same time step.  Thus, a fair computational comparison of various methods requires the dual setup, where a different time step size is used for each of them, so that all yield the same level of numerical error.  The time step 
that needs to be used for each method
can be 
estimated by means of the proposed tool.
\item Implementation of the proposed tool requires only 
the associated matrix pencils and thus, it allows quickly testing -- potentially new -- numerical schemes and interfacing strategies whose full implementation in the time-domain routine may be an involved procedure.  Schemes that do not fulfill the user’s numerical deformation
requirements 
can be discarded without further considerations.
\item The proposed technique is based on \ac{sssa} and thus, being strict, its results are valid around steady-state solutions $(\bfg x_o,\bfg y_o)$ of the system. Around $(\bfg x_o,\bfg y_o)$, comparison of $\bfb P(s)$, ${\bfb {P}}_{\scsym PC}(z)$ and ${\bfb {\bar P}}_{\scsym PC}(z)$ allows determining precisely the numerical deformation introduced by an integration method given the time step.
That said, the structure of the dynamic modes and stiffness of power system models, as well as the properties of integration methods, 
tend to be “robust” and therefore, 
our results provide a rough yet accurate quantification of the system's dynamics numerical deformation, also for changing operating conditions.  Similar considerations can be found in the existing literature, e.g.,~in \cite{arriaga:82_2,book:chow:13,tzounas2022tdistab,tzounas2022tdistabd}.
\end{itemize}

\subsection{Dense Matrix Formulation}
\label{sec:dense}
 
The matrix pencils ${\bfb {P}}_{\scsym PC}(z)$ and ${\bfb {\bar P}}_{\scsym PC}(z)$ were derived above using a sparse matrix formulation.  As discussed in Section~\ref{sec:dae:sssa}, this is the preferred option if a large-scale system is to be analyzed.  On the other hand, in Section~\ref{sec:dae:sssa} we also discussed that, for small/medium size systems, it is more efficient to work with dense matrices.  Thus, for the sake of completeness, we also provide an alternative way to compute the finite eigenvalues of the two pencils, based on dense matrices. 

From \eqref{eq:pc:g:lin}, $\tilde{\bfg y}_{n+1} = -\jac{g}{y}^{-1}\jac{g}{x}\tilde{\bfg x}_{n+1}$ and $\tilde{\bfg y}_n = -\jac{g}{y}^{-1}\jac{g}{x}\tilde{\bfg x}_n$, under the assumption that $\jac{g}{y}$ is non-singular.  Using these expressions, as well as \eqref{eq:corr:xt:lin}, in \eqref{eq:mem:corr:lin} (see Appendix), we find:
\begin{align}
 \label{eq:mem:pencil:dense:ext}
{\bfb P}_{\scsym PC}(z) &= z \bfb I  -
 (\bfb I + h \bfb C_r \AS )  \, , \\
  \label{eq:mem:pencil:dense:rig}
 {\bfb {\bar P}}_{\scsym PC}(z) &=
 z (\bfb I + \bfb M) -
(\bfb I + h \bfb C_r \AS + \bfb M ) \, ,
\end{align}
where $\bfb M = \frac{h}{2}\bfb C_{r-1}\jac{f}{y} \jac{g}{y}^{-1}\jac{g}{x}$ and $\bfb C_{r}$ is given by \eqref{eq:mem:cr}.

\subsection{Generic Adams-Bashforth PC Methods}

The proposed technique was described above for \ac{hm} but, in principle, is applicable to any \ac{psa}-based integration scheme.  
Since our focus is on \ac{pc} methods, we discuss here the application to a generic Adams-Bashforth \ac{pc} method.  An Adams-Bashforth method applied to \eqref{eq:dae} reads as:
%
\begin{equation}
\begin{aligned}
\label{eq:adams}
 \bfg \xpc^{(0)}_{n+1} &= \bfg x_{n}
 \hspace{-0.5mm}
   + h \sum_{j=0}^{k-1}\gamma_j \nabla^j \bfg f(\bfg x_{n},\bfg y_{n}) \, , \\ 
 \bfg \xpc^{(i)}_{n+1} &= \bfg x_{n}
 \hspace{-0.5mm}
 + h b_k \bfg f (\bfg \xpc^{(i-1)}_{n+1}, \bfg y_{\rm int}) 
\hspace{-0.5mm}
 + h 
 \hspace{-0.5mm}
 \sum_{j=0}^{k-1} b_j 
 \bfg{\mathcal{F}}_{n-k+j+1}
 \, , \\
\bfg x_{n+1} &= \bfg \xpc^{(r)}_{n+1} \, , \\
\bfg 0 &= h \bfg g(\bfg x_{n+1},\bfg y_{n+1}) \, .
  \end{aligned}    
\end{equation}
where $\nabla^j$ denotes the $j$-th order backward difference operator;  $\bfg{\mathcal{F}}_{n-k+j+1} = \bfg f(\bfg x_{n-k+j+1},\bfg y_{n-k+j+1})$;  and $k \in \mathbb{N}^*$.  Notice that \ac{hm} is a special case of the Adams-Bashforth method obtained from \eqref{eq:adams} for $k=1$, $\gamma_0=1$, and $b_0=b_1=0.5$.
%
Considering either perfect interfacing or extrapolation, the properties of \eqref{eq:adams} can be seen through a linear system of difference equations in the form:
\begin{equation}
\label{eq:difference}
 \bfg{\mathcal{E}} \bfb y_{n+1} = \bfg{\mathcal{A}} \bfb y_{n} \, ,
\end{equation}
and thus, similarly to previous sections, by studying a matrix pencil in the form $z\bfg{\mathcal{E}} - \bfg{\mathcal{A}}$.  The proof of \eqref{eq:difference} is provided in Section~\ref{sec:proof:diff} of the Appendix. 


\subsection{Delay-Based Numerical Stability Analysis in \cite{milanopsa}}
\label{sec:psafm}

To the best of our knowledge, the only attempt to study the numerical deformation caused by \ac{psa} to the dynamic modes of power systems is the technique described in \cite{milanopsa}.  
The main idea in \cite{milanopsa} is that the numerical effect of interfacing by extrapolation can be seen by studying the system that arises if all the algebraic variables that appear in the differential equations of \eqref{eq:dae} are taken from the previous time step.  Applying this approach changes \eqref{eq:dae} to the following time-delay system:
\begin{equation}
\begin{aligned}
\label{eq:ddae}
{\bfg x}'(t)  &= \bfg f( \bfg  x(t) , \bfg y (t-h) ) \, , \\
\bfg 0 &=  \bfg g( \bfg  x(t) , \bfg y (t) ) \, .
\end{aligned}    
\end{equation}
Then, linearization of \eqref{eq:ddae} gives:
\begin{equation}
  \begin{aligned}
    \label{eq:ddae:lin}
    \tilde {\bfg x}'(t)  &= 
    \jac{f}{x} \tilde {\bfg x}(t)  + \jac{f}{y} \tilde {\bfg y}(t-h) \, , \\
    \bfg 0 &=\jac{g}{x} \tilde {\bfg x}(t) + \jac{g}{y} \tilde {\bfg y}(t) \, , 
  \end{aligned}    
\end{equation}
and the matrix pencil of \eqref{eq:ddae:lin} is:
\begin{equation}\label{eq:ddae:pencil:sparse}
    \bfb P_d(s) = s \bfb E - \begin{bmatrix}
\jac{f}{x} & \bfg 0 \\
\jac{g}{x} & \jac{g}{y} \\
\end{bmatrix} -
\begin{bmatrix}
\bfg 0 & \jac{f}{y} \\
\bfg 0 & \bfg 0 \\
\end{bmatrix} e^{-sh} \, , 
\end{equation}
or, alternatively, using a dense matrix formulation,
\begin{equation}\label{eq:ddae:pencil:dense}
\bfb P_d(s) = s \bfb I - \jac{f}{x} +
\jac{f}{y}\jac{g}{y}^{-1}\jac{g}{x} e^{-sh} \, .
\end{equation}
\begin{figure}[ht!]\centering
\resizebox{\linewidth}{!}{\includegraphics{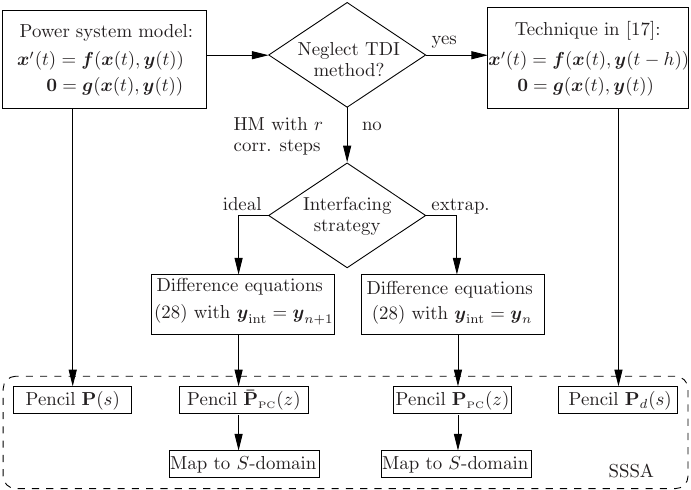}}
\caption{Procedure to obtain the 
matrix pencils considered.}
\label{fig:flow}
\end{figure}

\begin{table*}[!ht]
  \renewcommand{\arraystretch}{2.1}
  \renewcommand\tabcolsep{5.3pt}
  \centering
  \caption{Summary of matrix pencils considered. For each pencil, both the dense and sparse matrix form is given.  Numerical deformation can be determined by comparing the eigenvalues of ${\bfb P}_{\scsym PC}(z)$, $ {\bfb {\bar P}}_{\scsym PC}(z) $, {$\bfb P_d(s) $}, with those of ${\bfg P}(s)$. To facilitate comparison, eigenvalues of $Z$-domain pencils are mapped to the $S$-domain.}
  \label{tab:pencils}
  \begin{threeparttable}
    \begin{tabular}{lllll}
      \toprule\toprule
     {Approximation by} 
     & Symbol
     &
     {Sparse Matrix Form} 
     &
     Dense Matrix Form
    & 
    Domain~(Map to $S$)
       \\ 
      \midrule
  N/A~(DAE system) 
  & \centering {${\bfg P}(s)$}
  &  $s \bfb E - \bfb A$  
  &  $s \bfb I - \AS $
      & 
      $S$~(Not needed)
      \\
     \midrule
  \ac{hm}, 
  ${\bfg y}_{\rm int} = \bfg y_{n} \ \ $ 
  & {${\bfb P}_{\scsym PC}(z)$} 
  & 
  $z 
\renewcommand*{\arraystretch}{1.5}
\begin{bmatrix}
\bfb I & \bfg 0\\
\jac{g}{x}& \jac{g}{y} \\
\end{bmatrix} 
-  \begin{bmatrix}
\bfb I + h \bfb C_r\jac{f}{x} & h \bfb C_r\jac{f}{y}  \\
\bfg 0  &  \bfb 0 \\
\end{bmatrix}$ 
  & 
  $z \bfb I - (\bfb I + h \bfb C_r \AS )$ 
      & 
      $Z$ ($s=\frac{{\rm log}(z)}{h}$)
      \\
      \midrule
  \ac{hm}, ${\bfg y}_{\rm int} = \bfg y_{n+1} \ $ 
  & $ {\bfb {\bar P}}_{\scsym PC}(z) $
  & 
  $
  z
  \hspace{-1mm}
  \renewcommand*{\arraystretch}{1.5}
\begin{bmatrix}
\bfb I & 
-\bfb B \jac{f}{y}\\
\jac{g}{x} & \jac{g}{y} \\
\end{bmatrix}  
\hspace{-0.5mm}
- 
\begin{bmatrix}
\bfb I + h \bfb C_r \jac{f}{x} & 
\hspace{-1.2mm}
( h \bfb C_{r} 
\hspace{-0.7mm}
- \bfb B)\jac{f}{y} \\
\bfb 0 & \bfb 0 \\
\end{bmatrix}  
$
  & $ z (\bfb I + \bfb M) -
(\bfb I + h \bfb C_r \AS + \bfb M )$
& 
$Z$ ($s=\frac{{\rm log}(z)}{h}$)
\\
\midrule
PSA analysis in \cite{milanopsa} 
& {$\bfb P_d(s) $}
&  
\renewcommand*{\arraystretch}{1.5}
$s \bfb E - 
\begin{bmatrix}
\jac{f}{x} & \bfg 0 \\
\jac{g}{x} & \jac{g}{y} \\
\end{bmatrix} -
\begin{bmatrix}
\bfg 0 & \jac{f}{y} \\
\bfg 0 & \bfg 0 \\
\end{bmatrix} e^{-sh}$   
& $s \bfb I - \jac{f}{x} +
\jac{f}{y}\jac{g}{y}^{-1}\jac{g}{x} e^{-sh}$
& 
$S$~(Not needed)
      \\
\bottomrule\bottomrule
    \end{tabular}
  \end{threeparttable}
\end{table*}

Therefore, in terms of matrix pencils, the idea in \cite{milanopsa} is that the numerical effect of \ac{psa} (in particular the interface error) can be seen by comparing the eigenvalues of $\bfb P_d(s)$ and $\bfb P(s)$.  A limitation of this approach is that it does not take into account the integration scheme applied, or, equivalently, it assumes that its effect can be neglected.  For example, using \ac{hm} with $r=1$ or $2$ would make no difference on the interface error in view of $\bfb P_d(s)$.  On the contrary, the matrix pencils proposed in this paper to quantify the interface error, i.e.~${\bfb {\bar P}}_{\scsym PC}(z)$ and ${\bfb {P}}_{\scsym PC}(z)$, take into account both the numerical method applied and the adopted interfacing strategy, and, from this viewpoint, they are \textit{exact}.  

Finally, a summary of the pencils $\bfb P(s)$, ${\bfb {P}}_{\scsym PC}(z)$, ${\bfb {\bar P}}_{\scsym PC}(z)$, and $\bfb P_d(s)$, is given in Table~\ref{tab:pencils}. Moreover, a flow chart that outlines the process to derive these pencils is presented in Fig.~\ref{fig:flow}.  A comparison of the information given by these pencils is provided in the case study presented in Section~\ref{sec:case}.

\section{Case Studies}
\label{sec:case}

In this section, we illustrate the proposed technique to study the numerical stability and accuracy of \ac{psa}.  The results in Section~\ref{sec:39bus} are based on the 
IEEE 39-bus benchmark system, whereas 
Section~\ref{sec:aiits} is based on a 1479-bus model of the \acf{aiits}.  Simulations are carried out using the power system analysis software tool Dome \cite{vancouver}.  The examined systems permit efficient calculation of their whole spectra, and thus, throughout this section, eigenvalues
are computed with LAPACK \cite{lapack}, using dense matrix forms (see Table~\ref{tab:pencils}). 

\subsection{IEEE 39-bus system}
\label{sec:39bus}

This section considers the IEEE 39-bus 
system, the static and dynamic data of which are detailed in \cite{web:39bus}.  The 39-bus system comprises 10 synchronous machines represented by fourth order, two-axis models, 34 transmission lines, 12 transformers, and 19 loads.  Each machine is equipped with an \ac{avr}, a \ac{tg}, and a \ac{pss}. 
The \ac{dae} model has in total 129 states and 262 algebraic variables.  An equilibrium of the \acp{dae} is obtained from the solution of the power flow problem and the initialization of dynamic components.  Then, the small-signal dynamics 
of the system are determined by computing the eigenvalues of the matrix pencil $\bfb P(s)$ (see Section~\ref{sec:dae:sssa}).  In particular, eigenvalue analysis of $\bfb P(s)$ shows that the test system is stable around the examined equilibrium. 
 The system's fastest and slowest dynamics are represented, respectively, by the eigenvalues
$-106.01$ and $-0.02$, which gives from
\eqref{eq:stiffness}
a stiffness ratio of $5.3 \cdot 10^3$. 

\begin{figure*}[ht!]
  \centering
  \begin{subfigure}{.3\linewidth}
    \centering
    \resizebox{\linewidth}{!}{\includegraphics{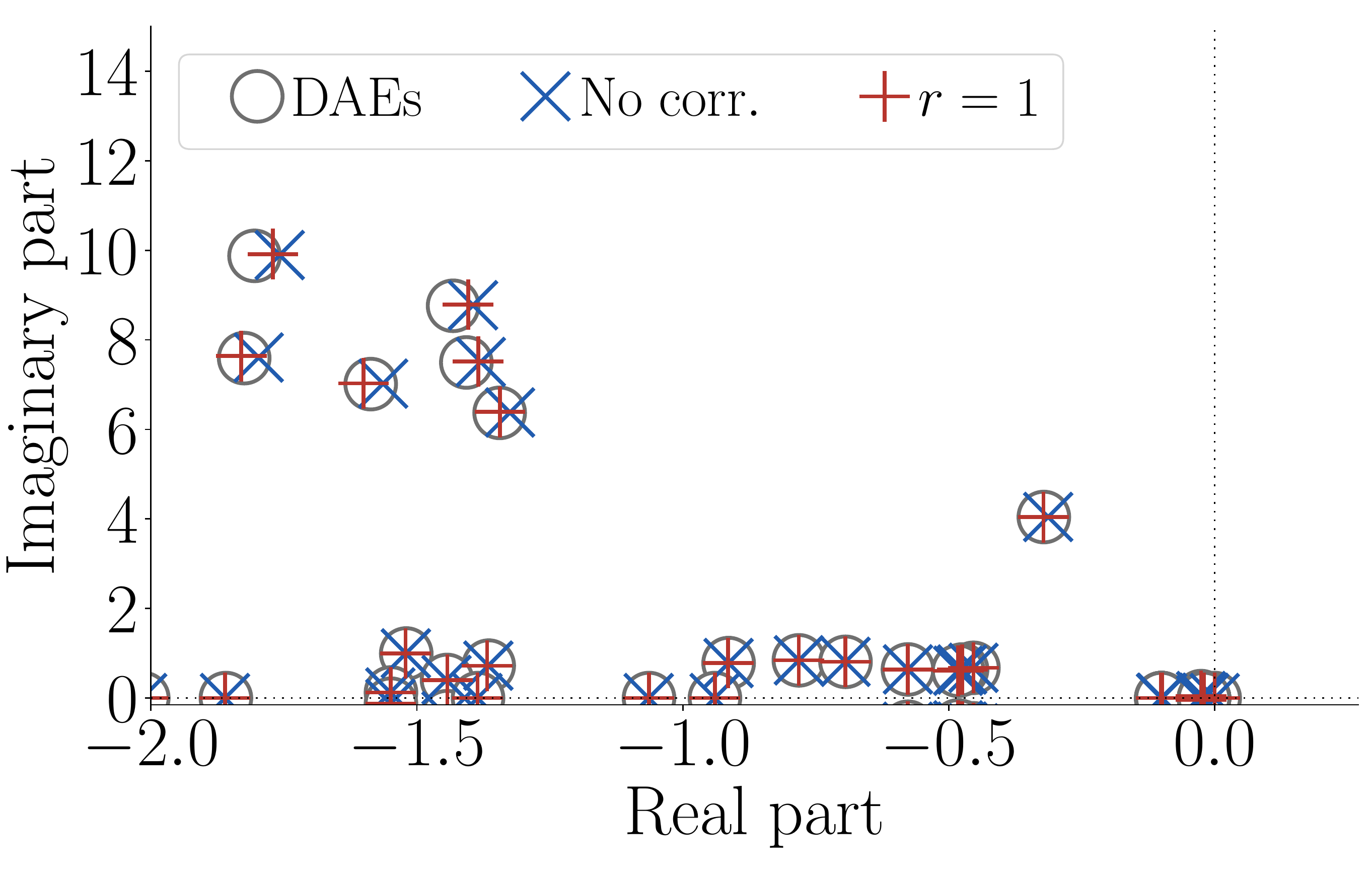}}
    \caption{$h=0.001$~s.}
    \label{fig:mem:ext:p001}
  \end{subfigure}
  \begin{subfigure}{.3\linewidth}
    \centering
    \resizebox{\linewidth}{!}{\includegraphics{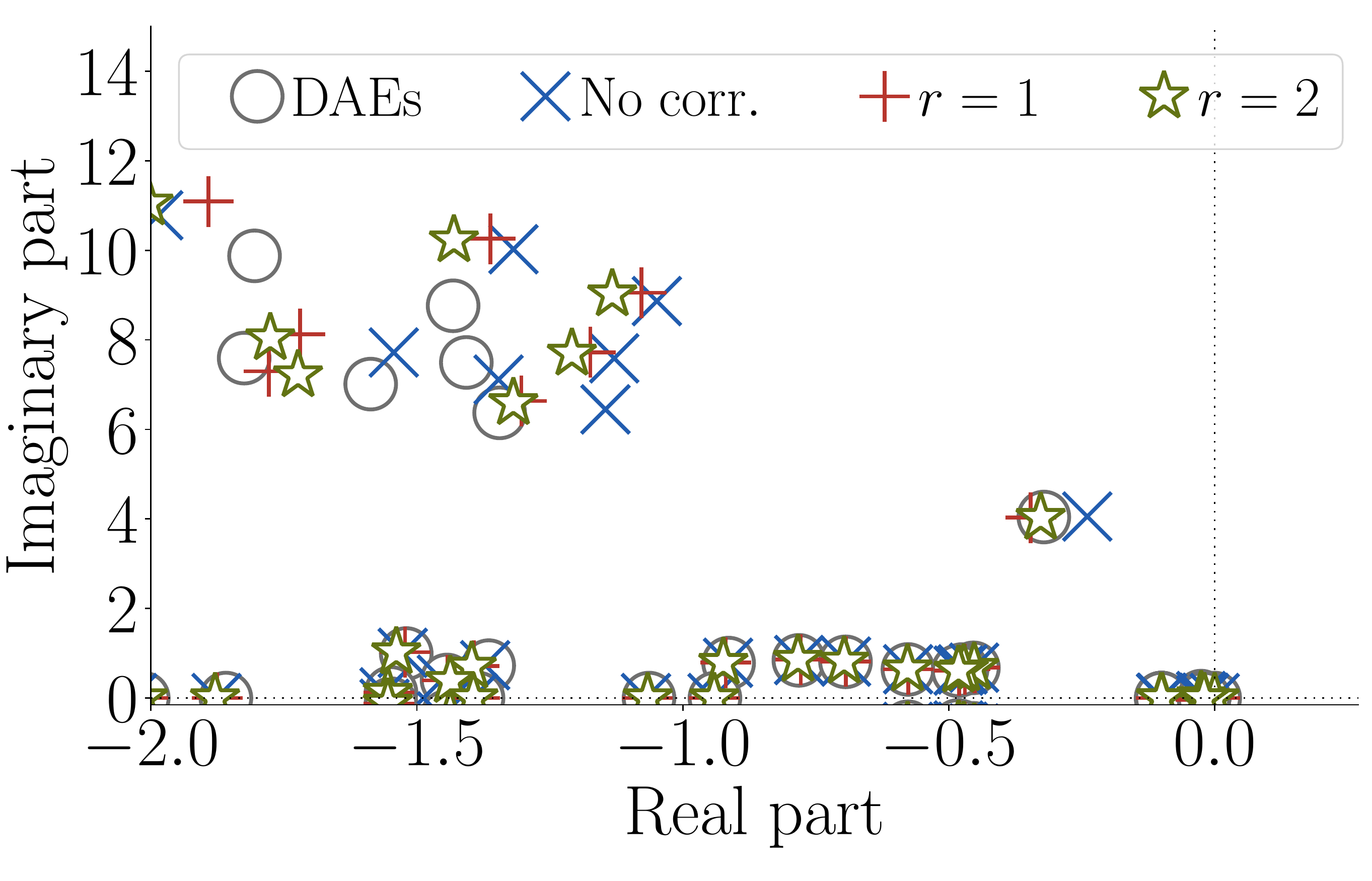}}
    \caption{$h=0.01$~s.}
    \label{fig:mem:ext:p01}
  \end{subfigure}
  \begin{subfigure}{.3\linewidth}
    \centering
    \resizebox{\linewidth}{!}{\includegraphics{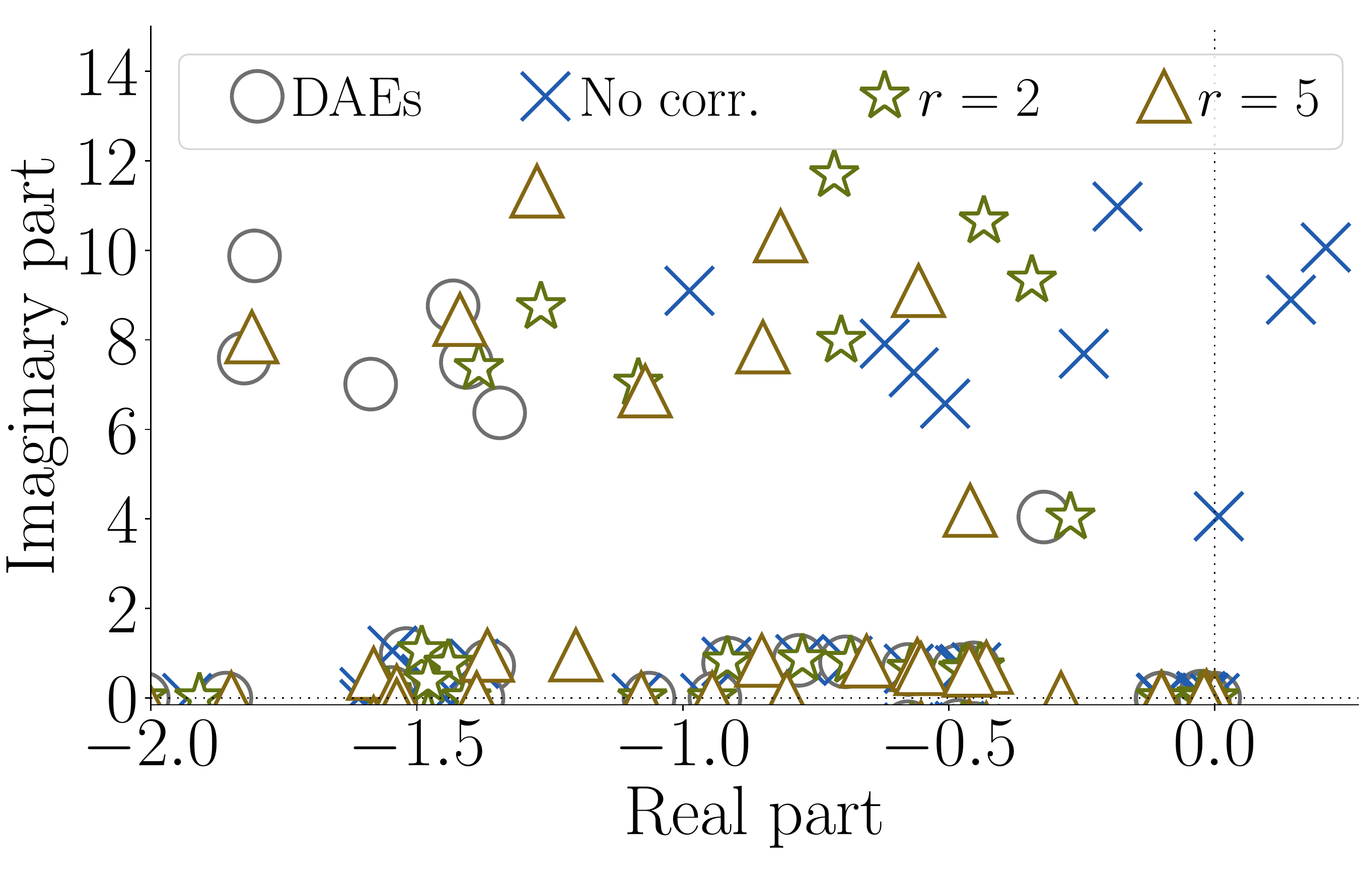}}
    \caption{$h=0.04$~s.}
    \label{fig:mem:ext:p04}
  \end{subfigure}
  \caption{39-bus system: Eigenvalue analysis of \ac{hm}, interfacing by extrapolation.}
  \label{fig:mem:ext}
  \vspace{-3mm}
\end{figure*}
\begin{figure*}[ht!]
  \centering
  \begin{subfigure}{.3\linewidth}
    \centering
    \resizebox{\linewidth}{!}{\includegraphics{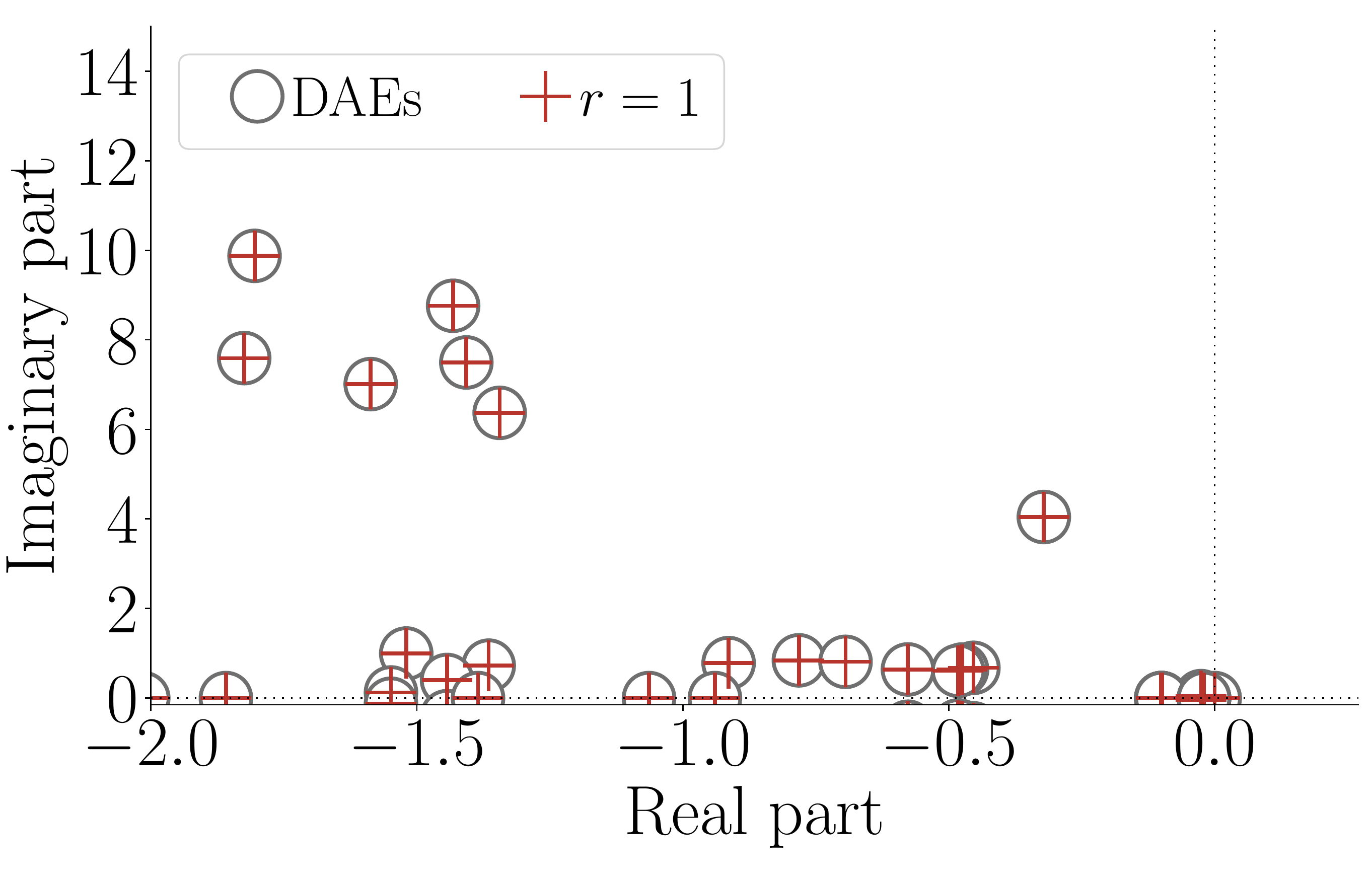}}
    \caption{$h=0.001$~s.}
    \label{fig:mem:pin:p001}
  \end{subfigure}
  \begin{subfigure}{.3\linewidth}
    \centering
    \resizebox{\linewidth}{!}{\includegraphics{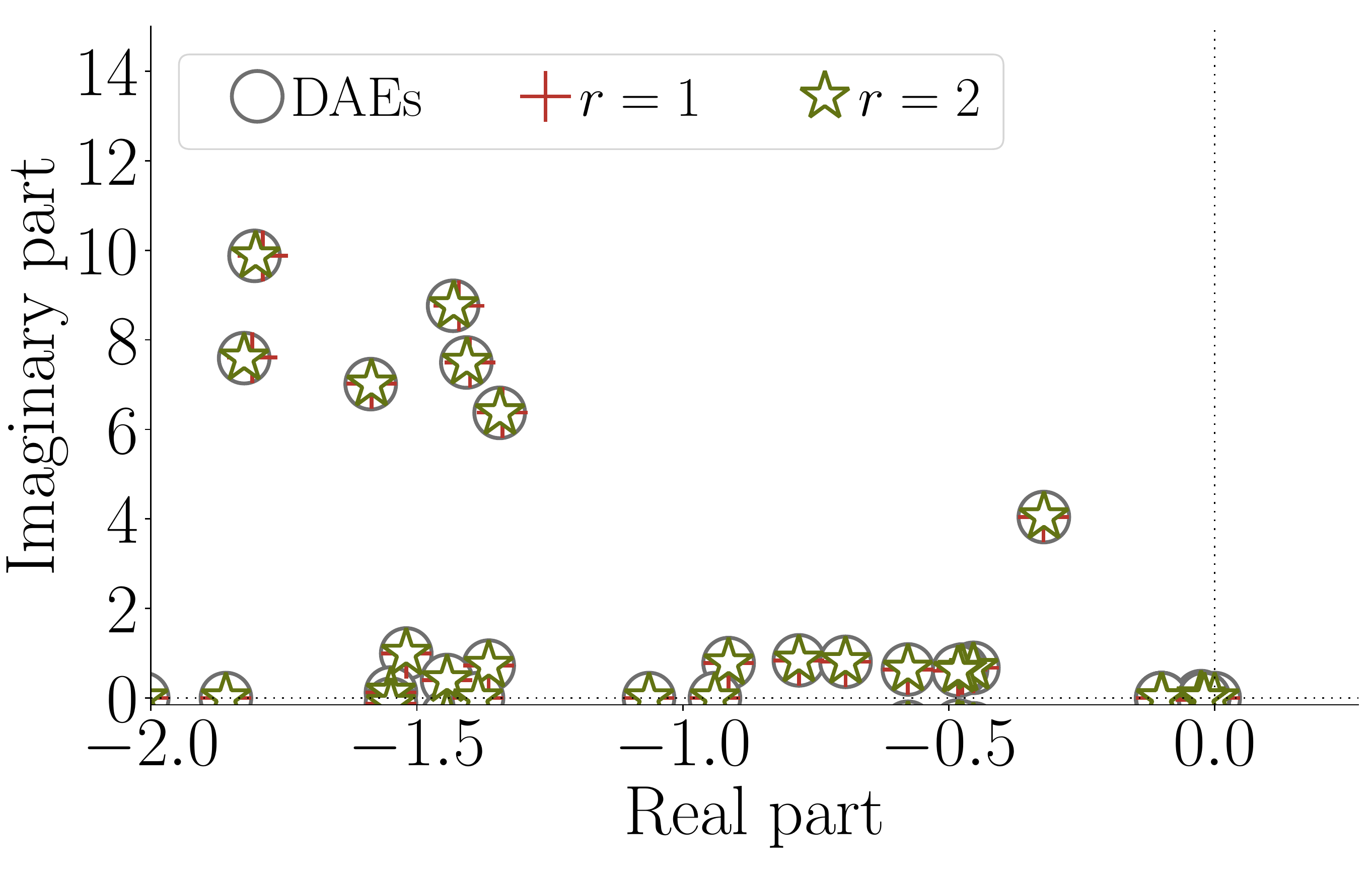}}
    \caption{$h=0.01$~s.}
    \label{fig:mem:pin:p01}
  \end{subfigure}
  \begin{subfigure}{.3\linewidth}
    \centering
    \resizebox{\linewidth}{!}{\includegraphics{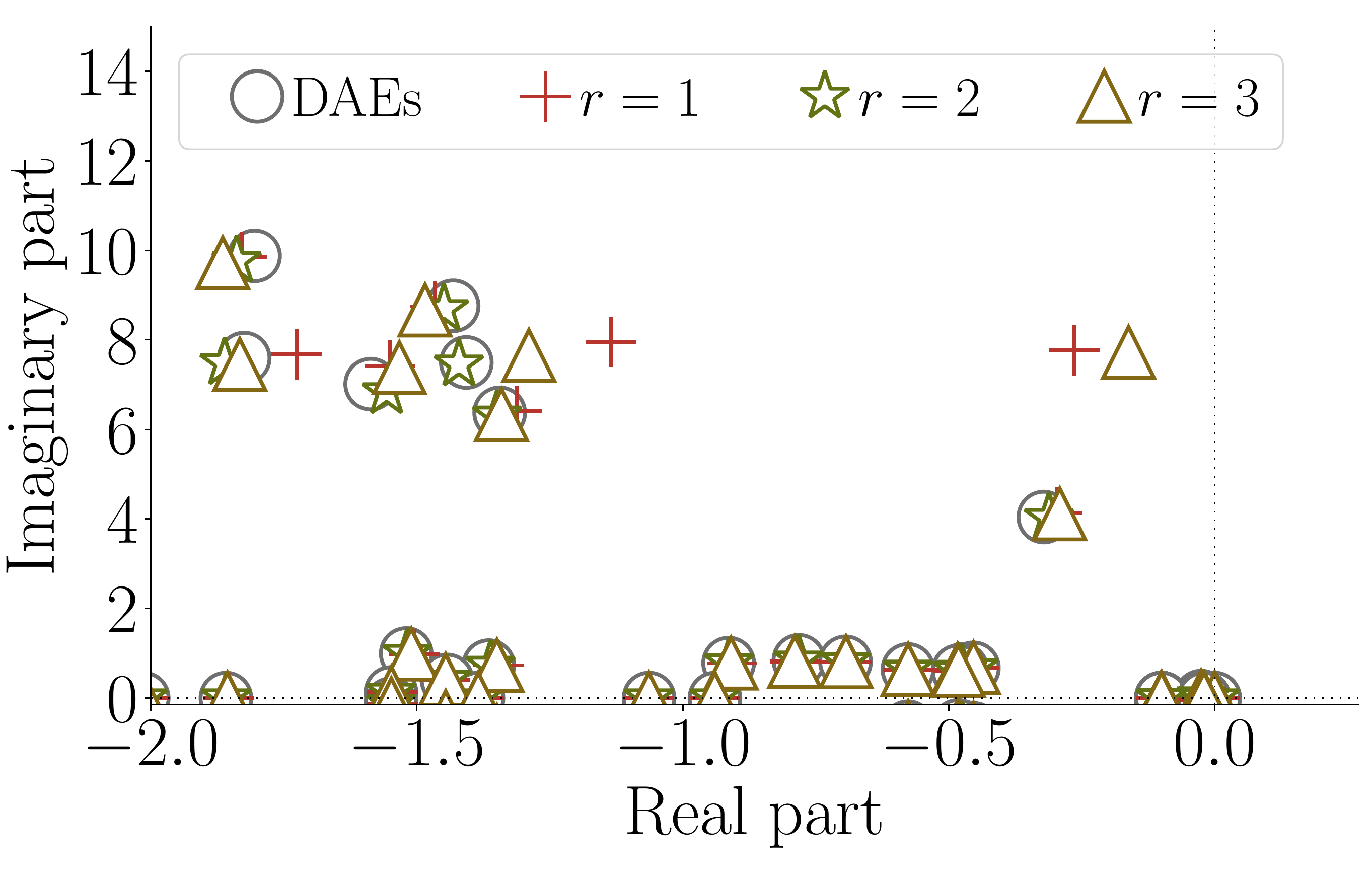}}
    \caption{$h=0.04$~s.}
    \label{fig:mem:pin:p04}
  \end{subfigure}
  \caption{39-bus system: Eigenvalue analysis of \ac{hm}, 
  perfect interfacing.}
  \label{fig:mem:pin}
  \vspace{-3mm}
\end{figure*}

\subsubsection{Interfacing by Extrapolation}
\label{sec:case:extr}

In the following, we consider that \ac{psa} is employed for the numerical integration of the system. In particular, we consider that the system is numerically solved using \ac{hm}, and that interfacing is handled through extrapolation of the known values of the algebraic variables, i.e.~${\bfg y}_{\rm int} = \bfg y_{n}$, as described in Section~\ref{sec:psa}.  To quantify the deformation caused to the dynamic modes of the system due to the application of the \ac{psa}-\ac{hm} numerical scheme, we calculate the pencil ${\bfb P}_{\scsym PC}(z)$ given by \eqref{eq:mem:pencil:dense:ext} and compute its eigenvalues, which we then map from the $Z$-domain to the $S$-domain according to \eqref{eq:z2s}.  Figure~\ref{fig:mem:ext} shows the numerical deformation introduced to the rightmost eigenvalues of the system for different time step sizes and for different number of corrector steps.  For the sake of comparison, the case where no corrector step is applied is also provided in each plot. Note that in this case the integration scheme reduces to \ac{fem}, see \eqref{eq:mem}.  Results suggest that if a small time step is employed, then a single corrector step ($r=1$) suffices to refine accuracy, see Fig.~\ref{fig:mem:ext:p001}, where $h=0.001$~s.  For $h=0.01$~s, one corrector step is not adequate to achieve precise representation of the system's response, while using $r=2$ slightly improves the accuracy for most modes, yet leaving a non-negligible error.  Moreover, higher values of $r$ increase the computational burden of the numerical method, without however providing a proportionate improvement in precision.  This is to be expected for large-enough time steps, where the effect of the extrapolation of algebraic variables becomes significant. Finally, for $h=0.04$~s, the solution is guaranteed to diverge if no corrector step is applied, since this case shows distorted eigenvalues in the right half of the $S$-plane, indicating numerical instability.  In particular, the numerical stability margin of the system when no corrector is applied is found to be $h=0.035$~s.  Including corrector steps, although it prevents numerical instability, does not yield a precise approximation of the system's dynamic modes.  Furthermore, interestingly, there also exist modes for which the accuracy deteriorates with higher values of $r$.

\begin{figure}[ht!]
\begin{subfigure}{1\linewidth}
    \centering
    \resizebox{0.9\linewidth}{!}{\includegraphics{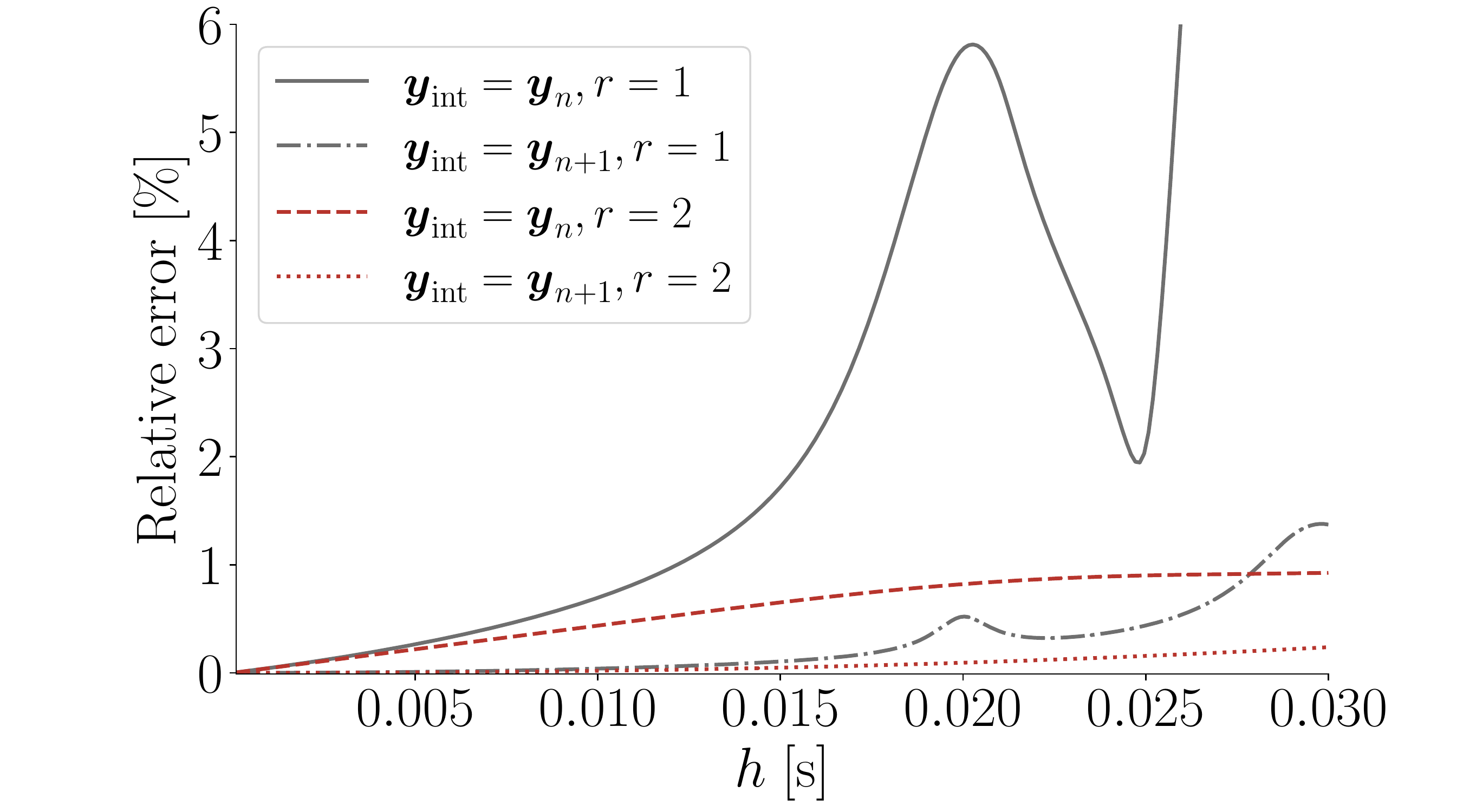}}
    \caption{Mode~1.}
    \label{fig:mem_m1}
  \end{subfigure}
  \begin{subfigure}{1\linewidth}
    \centering
    \resizebox{0.9\linewidth}{!}{\includegraphics{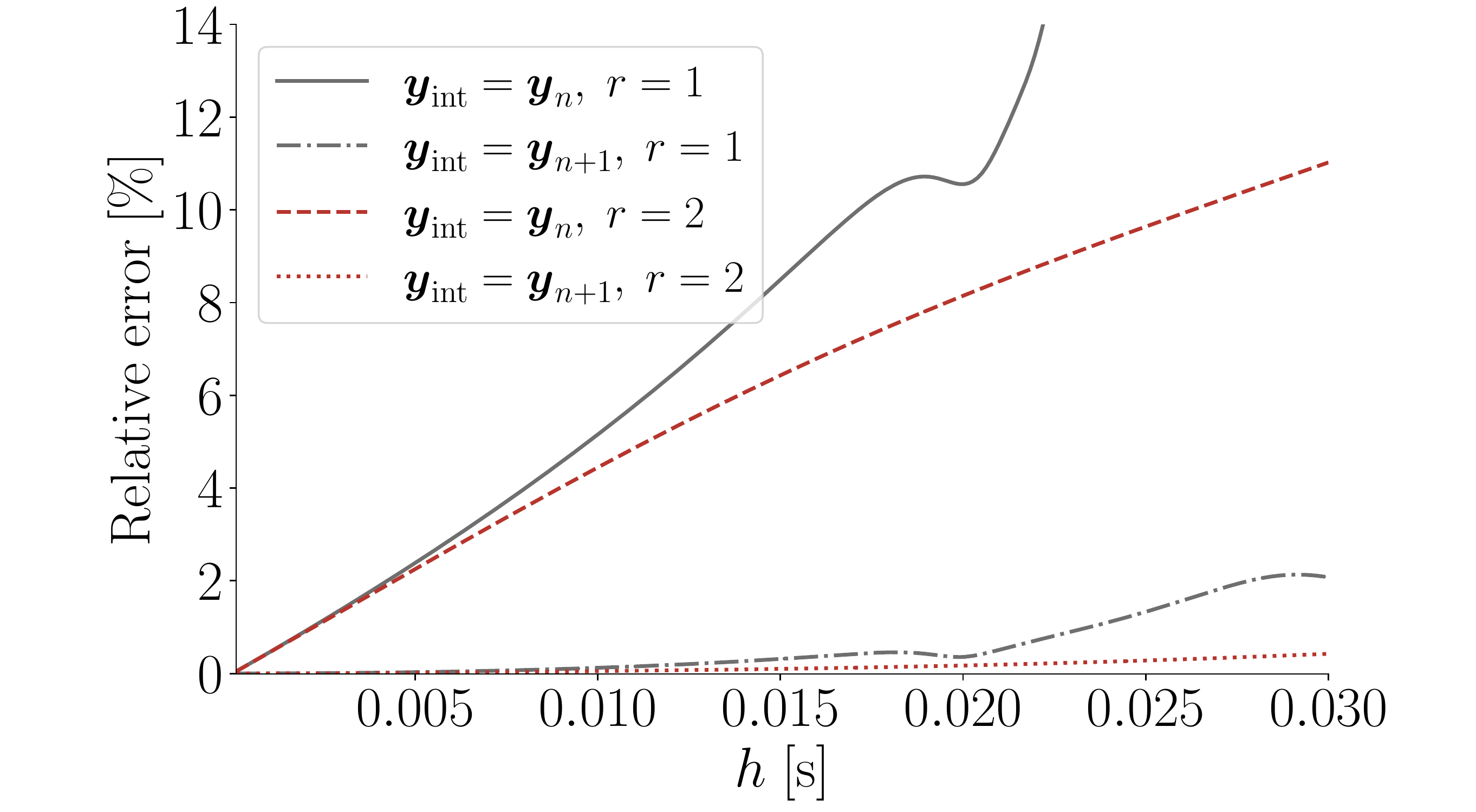}}
    \caption{Mode~2.}
    \label{fig:mem_m2}
  \end{subfigure}
\caption{9-bus system: Numerical deformation of least damped modes by \ac{hm}.}
  \label{fig:mem_m1m2}
  \vspace{-3mm}
\end{figure}

\subsubsection{Comparison with Perfect Interfacing}

In this section, we assume that the employed \ac{psa}-\ac{hm} scheme achieves elimination of the interface error or, equivalently, that ${\bfg y}_{\rm int} = \bfg y_{n+1}$.  The numerical deformation of the dynamic modes in this scenario can be quantified through the eigenvalues of the pencil ${\bfb {\bar P}}_{\scsym PC}(z)$.  The results for different time steps and number of corrector steps 
are presented in Fig.~\ref{fig:mem:pin}. These results indicate that, for time steps in the range $0.001$-$0.01$~s, using 1 to 2 corrector steps is adequate to achieve a precise approximation of the examined system's dynamics (Figs.~\ref{fig:mem:pin:p001},
\ref{fig:mem:pin:p01}).  For larger step sizes (see Fig.~\ref{fig:mem:pin:p04}), corrector steps allow only partial accuracy refinement.

In the following, we provide a further comparison of the results obtained with ${\bfg y}_{\rm int} = \bfg y_{n}$ and ${\bfg y}_{\rm int} = \bfg y_{n+1}$.  
With this aim, we focus on the two most critical modes of the system as defined by their damping ratios.  The most lightly damped mode (hereafter Mode~1) of the system is represented by the complex pair $-0.3212\pm \jj 4.0435$, which has damping ratio $7.92$\% and natural frequency $0.64$~Hz.  The second most lightly damped mode (hereafter Mode~2) is represented by the pair $-1.4318 \pm\jj 8.7610$, which has damping ratio $16.13$\% and natural frequency $1.41$~Hz. 

\begin{table}[!ht]
  \renewcommand{\arraystretch}{1.2}
  \centering
  \caption{39-bus system: Upper time step bounds that satisfy a $0.1$\% relative error requirement for Modes~1 and 2.}
  \label{tab:39bus:hmax}
  \begin{threeparttable}
    \begin{tabular}{c|cccc}
      \toprule\toprule
    $r$ & 1 & 1 & 2 & 2 \\
     \midrule
     $\bfg y_{\rm int}$ & $\bfg y_{n}$ & $\bfg y_{n+1}$ & $\bfg y_{n}$ & $\bfg y_{n+1}$\\
      \midrule
     max. time step~[s]
    & 0.0002 &  0.009 & 0.0002 & 0.015 \\
\bottomrule\bottomrule
    \end{tabular}
  \end{threeparttable}
\end{table}

We track the numerical deformation of Modes~1 and 2 as functions of the integration time step $h$.  The results for two different setups of the \ac{hm} ($r=1$ and $2$) are summarized in Fig.~\ref{fig:mem_m1m2}, where the relative error is defined as $100\cdot{|\hat s_i-s_i|}/{|s_i|}$, with $s_i$ being an eigenvalue of the system and $\hat s_i$ the corresponding numerically deformed eigenvalue.  Figure~\ref{fig:mem_m1} shows that, although extrapolating algebraic variables from the previous time step (${\bfg y}_{\rm int} = \bfg y_{n}$) introduces distortion, the structure of the integration method plays the most important role in the approximation of Mode~1.  In particular, when $r=1$, the interface error causes a deformation that becomes significant for large values of $h$.  Yet, with two corrector steps ($r=2$), the error is drastically reduced for all values of $h$ considered, with the relative deformation being below $1$\%, regardless of the adopted interfacing strategy.  On the other hand, Fig.~\ref{fig:mem_m2} indicates that the interfacing strategy plays the most important role in the approximation of Mode~2.  Specifically, extrapolation of algebraic variables significantly deforms Mode~2 both for small and larger time steps, with additional corrector steps not being able to notably reduce the error. 
Finally, assuming for the two modes a prescribed relative error requirement smaller than $0.1$\%, the estimated maximum admissible time step for the four different \ac{psa} setups examined are extracted 
from Fig.~\ref{fig:mem_m1m2} and summarized in
Table~\ref{tab:39bus:hmax}.  We note here that
the above discussion is based on two modes 
for the sake of 
illustration, but this is not a limitation of the proposed analysis, which can be extended to include more, or even all system modes.

\subsubsection{Comparison with Method in \cite{milanopsa}}

\begin{figure}[ht!]
  \begin{subfigure}{1\linewidth}\centering
    \resizebox{0.9\linewidth}{!}{\includegraphics{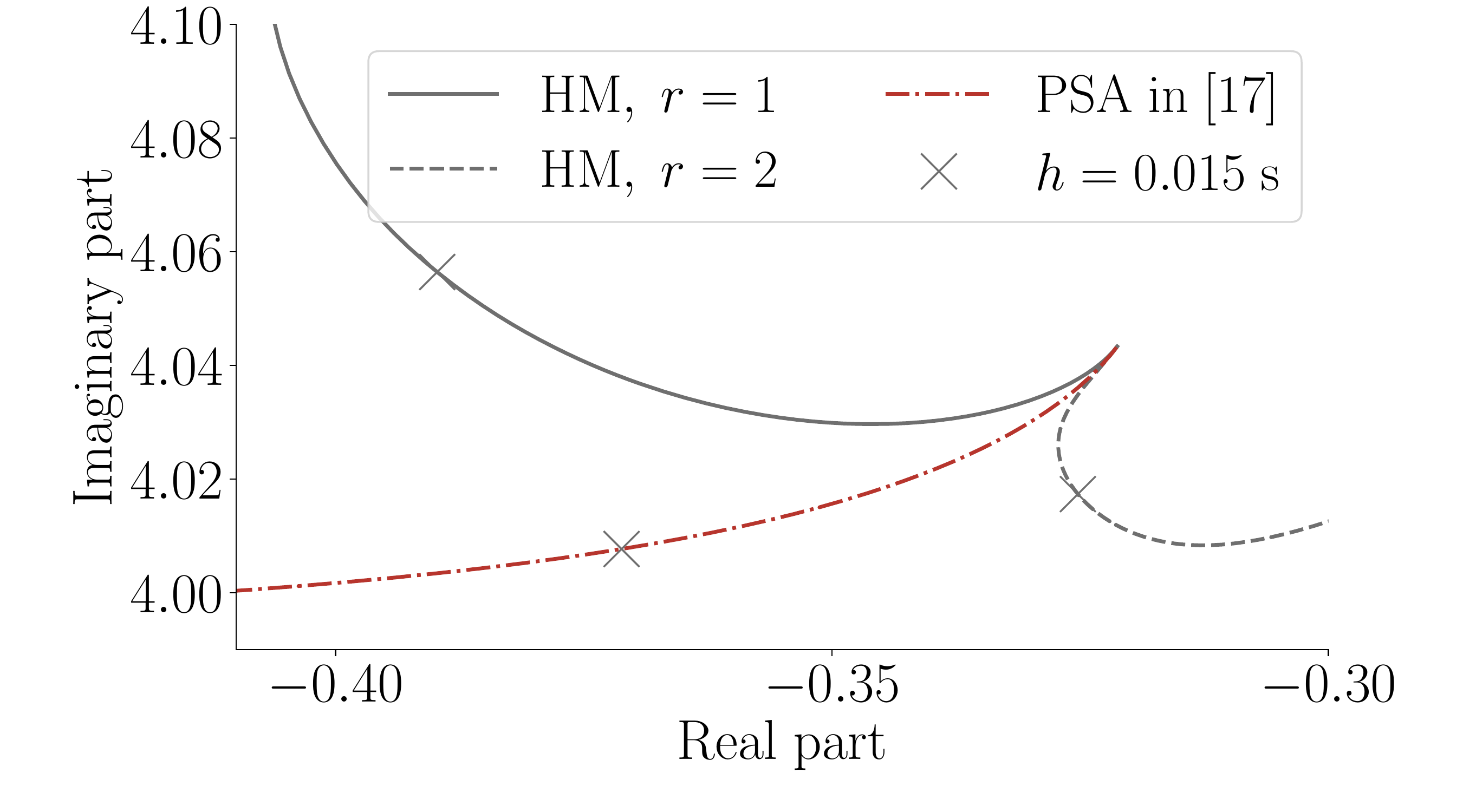}}
    \vspace{-2mm}
    \caption{Mode~1.}
    \label{fig:psa_m1}
  \end{subfigure}
  \begin{subfigure}{1\linewidth}\centering
    \resizebox{0.9\linewidth}{!}{\includegraphics{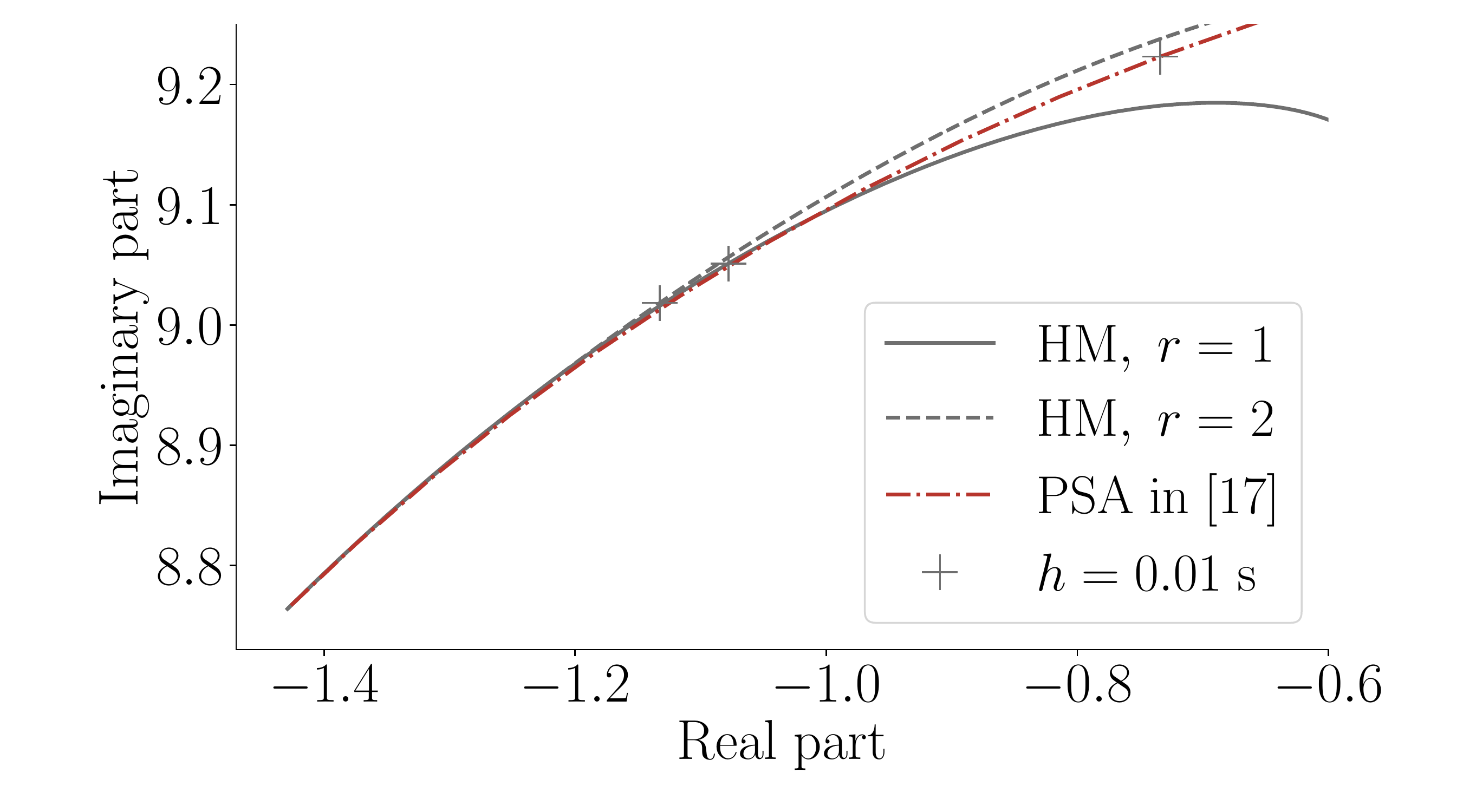}}
    \vspace{-2mm}
    \caption{Mode~2.}
    \label{fig:psa_m2}
  \end{subfigure}
\caption{39-bus system: HM~(${\bfg y}_{\rm int} = \bfg y_{n}$) vs method in \cite{milanopsa}.}
  \label{fig:psa:fm}
  \vspace{-3mm}
\end{figure}

We finally provide a discussion on the accuracy of the delay-based analysis of the \ac{psa} described in \cite{milanopsa}.  As mentioned in Section~\ref{sec:psafm}, \cite{milanopsa} focuses on the effect of extrapolating the algebraic variables while, on the other hand, effectively neglecting the impact of the specific integration method applied.  To illustrate how this technique compares to the approach described in this paper, we calculate the matrix pencil $\bfb P_d(s)$ and compute its eigenvalues as the simulation time step $h$ is varied.  Figure~\ref{fig:psa:fm} presents the root loci obtained for the two least damped modes of the system (Modes~1 and 2).  From the figure, it is obvious that the method in \cite{milanopsa} is not able to capture modes whose deformation is more due to the integration method applied and less due to the extrapolation of algebraic variables, such as Mode~1 (see Fig.~\ref{fig:psa_m1}).  On the contrary, it can be accurate in tracking the approximation of modes like Mode~2, the deformation of which is mostly due to interfacing by extrapolation.

\subsection{All-Island Irish Transmission System}
\label{sec:aiits}

This section is based on a 1479-bus model of the \ac{aiits}.  Topology and steady-state 
data have been provided by the Irish transmission system operator, EirGrid Group, whereas dynamic data have been determined based on our knowledge about the current technology of generators and automatic controllers.  The system model consists of 796~lines, 1055~transformers, 245~loads, 176~wind generators, and 22~\acp{sg} equipped with \acp{avr} and \acp{tg}.  Moreover, $6$~\acp{sg} are equipped with \acp{pss}.  In total, the system model includes 1615 state and 7225 algebraic variables.  Eigenvalue analysis of the matrix pencil $\bfb P(s)$ shows that the \ac{dae} power system model is stable under small disturbances.  The fastest and slowest eigenvalues of the system are respectively, $-99900.01$ and $-0.0013$.  These yield a stiffness ratio of $7.6 \cdot 10^7$.  Note that this is four orders of magnitude larger than the stiffness of the 39-bus system considered in Section~\ref{sec:39bus}.

\begin{figure}[ht!]
  \begin{subfigure}{1\linewidth}\centering
    \resizebox{0.92\linewidth}{!}{\includegraphics{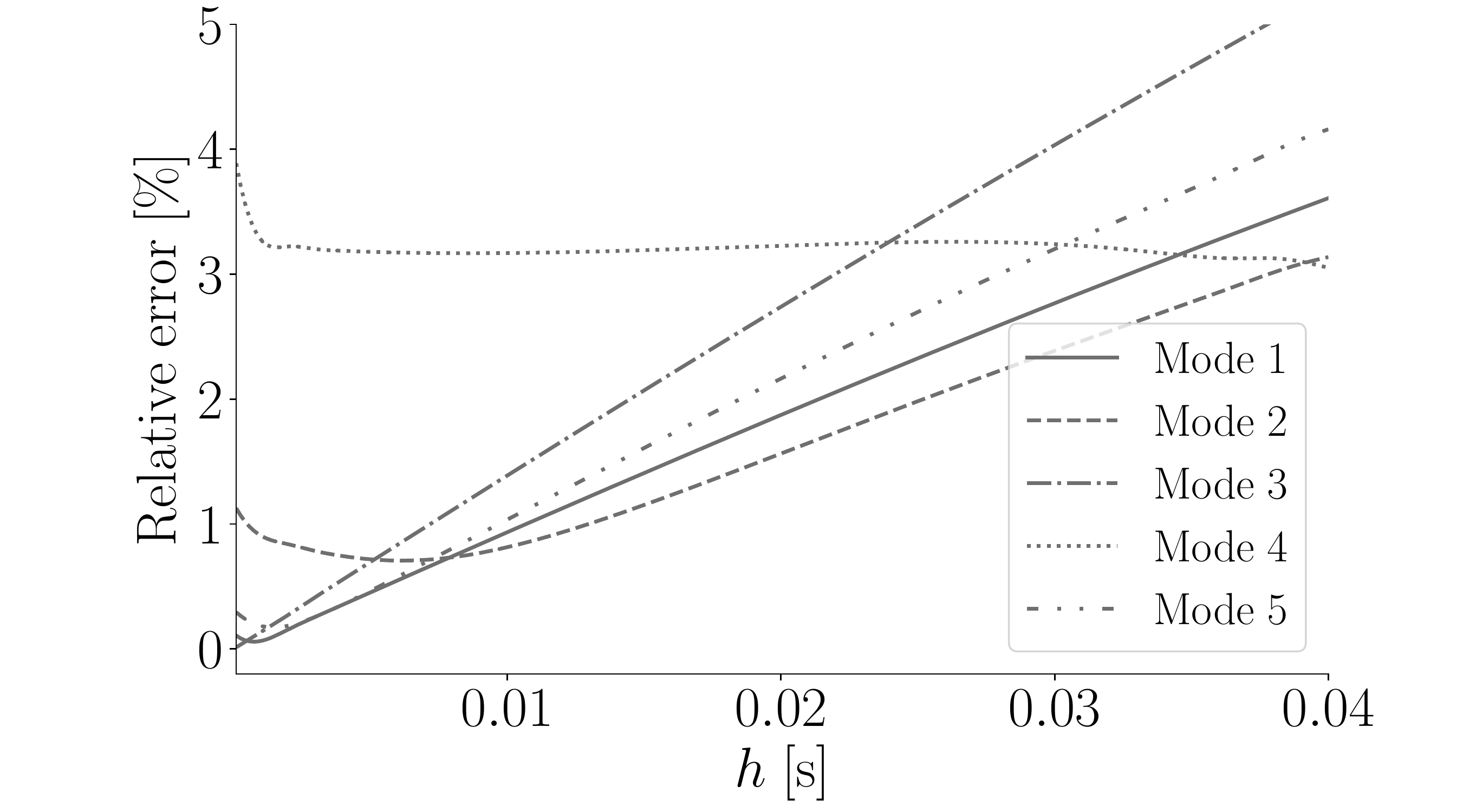}}
    \vspace{-1mm}
    \caption{Relative error.}
   \label{fig:aiits:ds}
  \end{subfigure}
  \begin{subfigure}{1\linewidth}\centering
    \resizebox{0.92\linewidth}{!}{\includegraphics{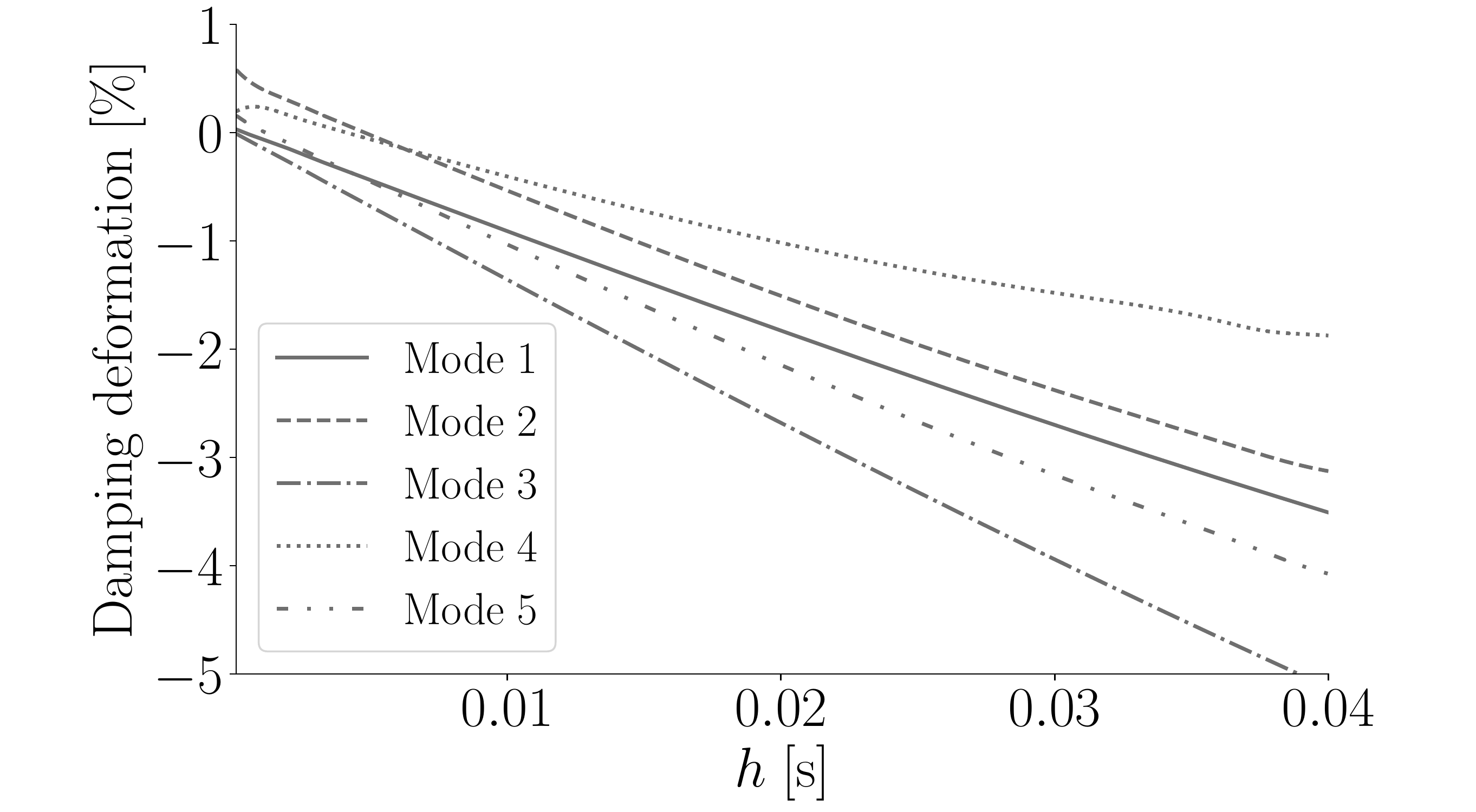}}
    \vspace{-1mm}
    \caption{Damping ratio deformation.}
    \label{fig:aiits:dz}
  \end{subfigure}
\caption{AIITS: Deformation of least damped modes by \ac{hm}~(${\bfg y}_{\rm int} = \bfg y_{n}$).}
  \label{fig:aiits}
  \vspace{-3mm}
\end{figure}

For the sake of illustration, we assume that \ac{hm} with interfacing by extrapolation of algebraic variables (${\bfg y}_{\rm int} = \bfg y_{n}$) is employed for the solution of the \ac{aiits}.  In this case, to apply the proposed technique we compute the eigenvalues of the associated pencil ${\bfb P}_{\scsym PC}(z)$ (see \eqref{eq:mem:pencil:dense:ext}). The full spectrum of ${\bfb P}_{\scsym PC}(z)$ for a given time step is 
computed with LAPACK in $\sim2.4$~s.\footnote{All simulations are carried out with a 64-bit Linux operating system running on a 8-core Intel i7 1.8~GHz, 8 GB laptop.}
We track the deformation of the five least damped dynamic modes of the system as the integration time step size is varied and we present the results in Figs.~\ref{fig:aiits}.  

\begin{figure}[ht!]
  \centering
  \includegraphics[width=0.49\textwidth]{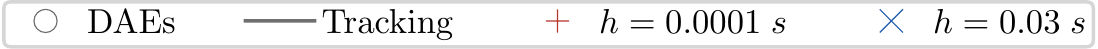}
  \begin{subfigure}{.49\linewidth}
    \hspace{-0.2cm}
    \centering
    \includegraphics[width=\textwidth]{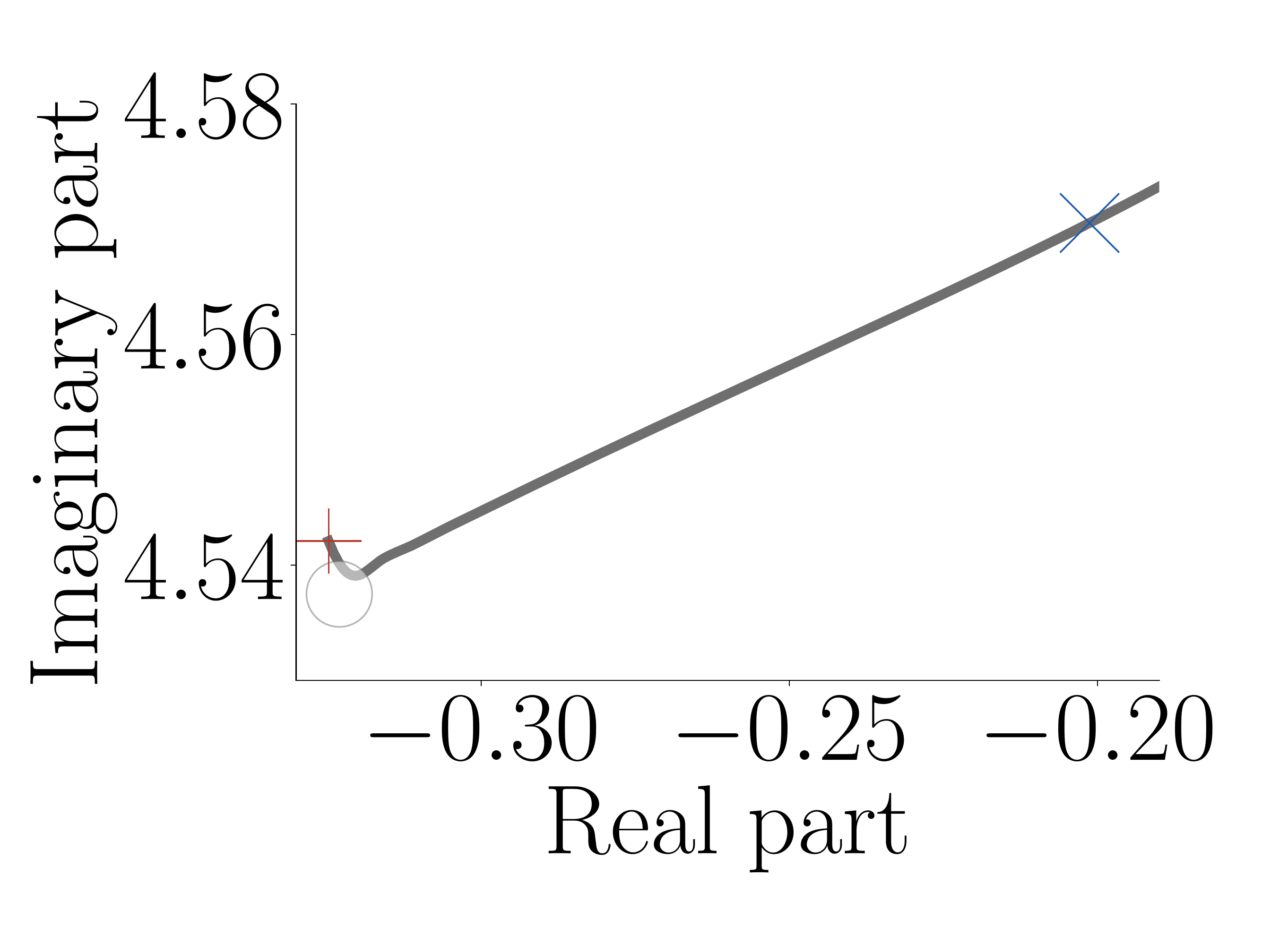}
    \caption{Mode 1 (${\bfg y}_{\rm int} = \bfg y_{n}$).}\label{fig:aiits:m1}
    \vspace{1mm}
  \end{subfigure}
  \begin{subfigure}{.49\linewidth}
    \centering
    \includegraphics[width=\textwidth]{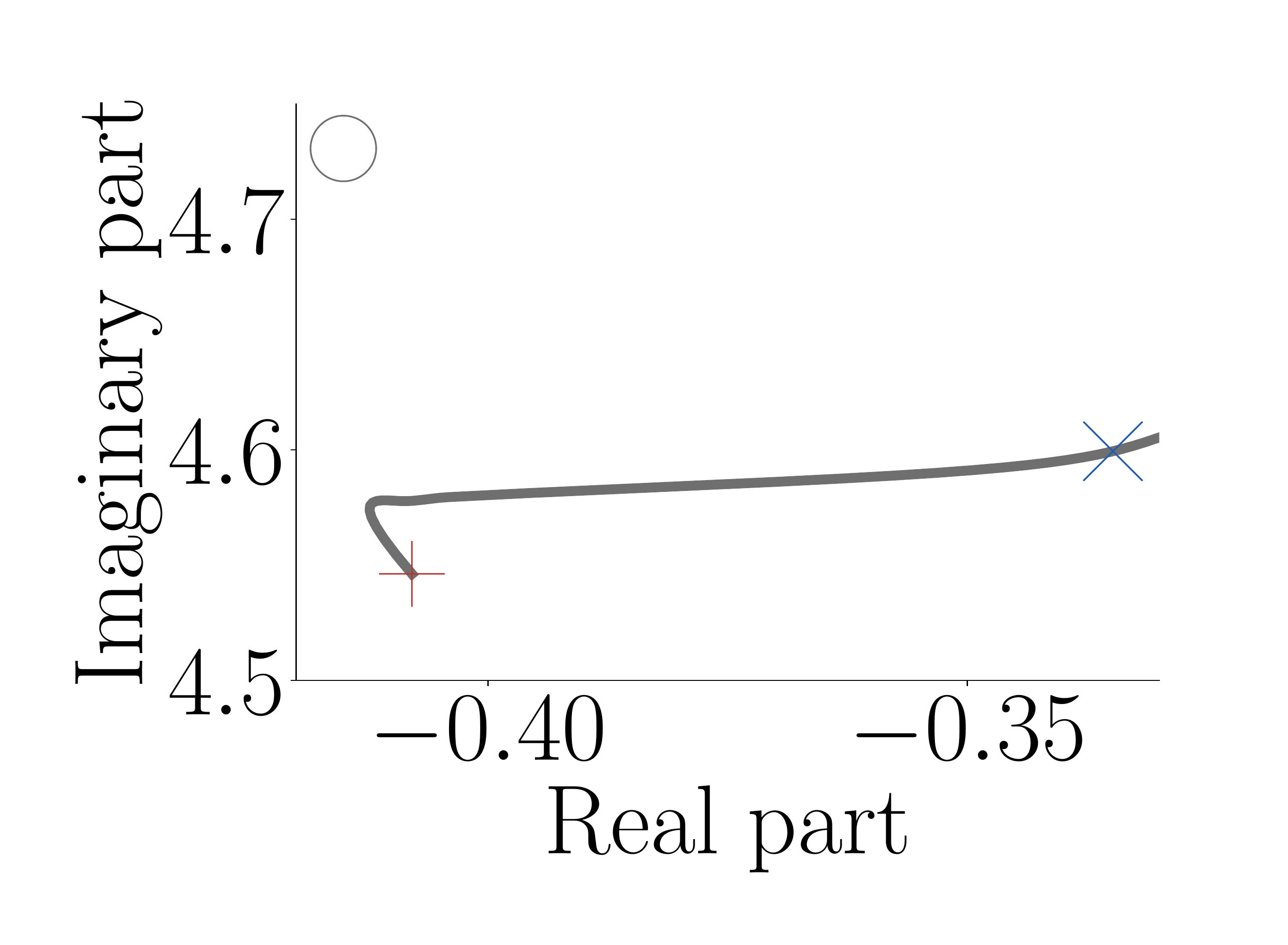}
    \caption{Mode 4 (${\bfg y}_{\rm int} = \bfg y_{n}$).}\label{fig:aiits:m4}
  \end{subfigure}
  \caption{AIITS: Eigenvalue tracking for Modes 1 and 4.}
  \label{fig:asbs_mem_aiits}
\end{figure}

The damping deformation shown in Fig.~\ref{fig:aiits:dz} is defined as follows.  If $s_i$ is an eigenvalue of the \ac{dae} system with damping ratio $\zeta$, and $\hat s_i$ is the corresponding deformed eigenvalue with damping ratio $\hat\zeta$, then the damping deformation is calculated as $100\cdot (\hat\zeta-\zeta)$.  Results indicate that numerical errors 
tend to increase with the time step, although they can also 
remain constant or decrease 
for some modes and in certain regions,
e.g.~see Modes~2 and 4 in 
Fig.~\ref{fig:aiits} 
Furthermore, we note that errors for some modes can be very large even for very small time steps (an
effect that was not observed on the IEEE 39-bus system).
For example, the relative error of Mode~4 for $h=10^{-4}$~s is about $4$\%.
This effect is caused by the one-step mismatch between state and algebraic variables present when interfacing by extrapolation.  This effect is further illustrated in Fig.~\ref{fig:asbs_mem_aiits}.


%

\section{Conclusion}
\label{sec:conclusion}

This paper proposes a novel technique based on \ac{sssa} to study the numerical stability and precision of the \ac{psa} applied for the solution of \ac{dae} power system models.  The proposed technique takes into account the dynamics of the model to be solved, the integration method applied, as well as the adopted interfacing strategy, and allows estimating useful time step bounds that achieve prescribed simulation accuracy criteria.  Formulation and simulation results are given for the well-known \ac{hm}, while the applicability of the proposed approach to any method of the Adams-Bashforth family is duly discussed.  Comparison with existing literature is also provided.  

We will dedicate future work to study the effect on the proposed tool of automatic time step size and order control techniques.
Moreover, we will employ the proposed approach to study
the numerical robustness of state-of-art power converter models 
solved with \ac{psa}, such as the ones in \cite{ramasubramanian2020positive,ramasubramanian2022parameterization}.

\appendix

\subsection{Proof of (\ref{eq:mem:predcor:lin})}
\label{sec:proof:mem}

Let first rewrite \eqref{eq:mem:pred:lin} and \eqref{eq:mem:cori:lin} as follows:
\begin{align}
 \label{eq:mem:pred:lin2}
\tilde {\bfg \xpc}^{(0)}_{n+1} &= 
(\bfb I +  h \jac{f}{x})
\tilde {\bfg x}_{n}
+ h \jac{f}{y} \tilde {\bfg y}_{n} \, , \\
\label{eq:mem:cori:lin2}
\tilde {\bfg \xpc}^{(i)}_{n+1} &= 
(\bfb I +  \frac{h}{2} \jac{f}{x})
\tilde {\bfg x}_{n} 
\hspace{-1mm}
+ \frac{h}{2} 
\jac{f}{x} \tilde {\bfg \xpc}^{(i-1)}_{n+1} 
\hspace{-1mm}
  + \frac{h}{2} 
  \jac{f}{y} (\tilde {\bfg y}_{n} + \tilde {\bfg y}_{\rm int}) \, .
\end{align}
For the $1$-st corrector step, using \eqref{eq:mem:pred:lin2} in \eqref{eq:mem:cori:lin2} for $i=1$ gives:
\begin{equation}
\begin{aligned}\label{eq:mem:cor1:lin:eq2}
\tilde {\bfg \xpc}^{(1)}_{n+1} = & 
( \bfb I + h \jac{f}{x}
+ 
\frac{h^2}{2}\jac{f}{x}^2 ) 
\tilde {\bfg x}_{n} \\ & 
+
 (\bfb I + h \jac{f}{x} ) \frac{h}{2} \jac{f}{y} \tilde {\bfg y}_{n} + \frac{h}{2} \jac{f}{y} \tilde {\bfg y}_{\rm int} 
 \, .
\end{aligned}\end{equation}
Similarly, for $i=2$, using \eqref{eq:mem:cor1:lin:eq2} in
\eqref{eq:mem:cori:lin2} gives:
%
%
%
\begin{equation}
\begin{aligned}
\nonumber
\tilde {\bfg \xpc}^{(2)}_{n+1} = & 
(\bfb I + (\bfb I 
+ \frac{h}{2} \jac{f}{x} +
\frac{h^2}{4} \jac{f}{x}^2) h \jac{f}{x})
\tilde {\bfg x}_{n} +
(\frac{1}{2} \bfb I + \frac{h}{4} \jac{f}{x} \\& + 
\frac{h^2}{4} \jac{f}{x}^2 ) 
{h} \jac{f}{y} \tilde {\bfg y}_{n} + 
(\frac{1}{2} \bfb I + \frac{h}{4} \jac{f}{x} ) 
{h} \jac{f}{y} \tilde {\bfg y}_{\rm int}  \, .
\end{aligned}\end{equation}
By induction, we get for the $r$-th corrector step: 
\begin{equation}
\begin{aligned}\label{eq:mem:corr:lin}
\tilde {\bfg \xpc}^{(r)}_{n+1} = 
(\bfb I + h \bfb C_r \jac{f}{x})
\tilde {\bfg x}_{n} 
&+ 
( \bfb C_{r} - 0.5 \bfb C_{r-1} )
h\jac{f}{y} \tilde {\bfg y}_{n}
\\ &+
0.5{h} \bfb C_{r-1}
\jac{f}{y} \tilde {\bfg y}_{\rm int} \, ,
\end{aligned}\end{equation}
where $\bfb C_{r}$ is given by \eqref{eq:mem:cr},
and 
$\bfb C_{r}= \bfb{I}+0.5h\bfb C_{r-1} \jac{f}{x}$.
From \eqref{eq:corr:xt:lin},~\eqref{eq:pc:g:lin},~\eqref{eq:mem:corr:lin}, we arrive at \eqref{eq:mem:predcor:lin}.

\subsection{Proof of (\ref{eq:difference})}
\label{sec:proof:diff}

Let us first rewrite the predictor in the following form:
\begin{align}
 \label{eq:adams:pred:2}
 \bfg \xpc^{(0)}_{n+1} &= \bfg x_{n}
  + h \sum_{j=0}^{k-1} c_j 
   \bfg f(\bfg x_{n-k+j+1},\bfg y_{n-k+j+1}) \, ,
\end{align}
where the coefficients $c_j$ can be obtained by using backward difference properties, i.e.:
\begin{align}
\nonumber
\nabla^{j+1} \bfg f(\bfg x_{n},\bfg y_{n}) &= \nabla^{j}\bfg f(\bfg x_{n},\bfg y_{n}) - \nabla^{j}\bfg f(\bfg x_{n-1},\bfg y_{n-1}) \, , \\
\nonumber
\nabla^0 \bfg f(\bfg x_{n},\bfg y_{n}) &= \bfg f(\bfg x_{n},\bfg y_{n})
 \, .
\end{align}
Additionally, the equilibrium point $(\bfg x_o,\bfg y_o)$ is also a fixed point of the numerical method, under the assumption that the system has been in steady-state for a time equal to $kh$.  Linearizing \eqref{eq:adams} around this point, we have for the predictor:
\begin{equation}
 \label{eq:adams:pred:lin}
\tilde {\bfg \xpc}^{(0)}_{n+1} = 
\tilde {\bfg x}_{n}
	+ h \sum_{j=0}^{k-1} c_j 
	(\jac{f}{x} \tilde {\bfg x}_{n-k+j+1} + \jac{f}{y} \tilde {\bfg y}_{n-k+j+1} ) 
	\, , 
\end{equation}
and for the $i$-th corrector step:
%
\begin{align}
 \nonumber
\tilde {\bfg \xpc}^{(i)}_{n+1} = 
\tilde {\bfg x}_{n}
 & + h \sum_{j=0}^{k-1} b_j (\jac{f}{x} \tilde {\bfg x}_{n-k+j+1} + \jac{f}{y} \tilde {\bfg y}_{n-k+j+1} ) 
\\
\label{eq:adams:corr:lin}
& + h b_{k} (\jac{f}{x} \tilde {\bfg \xpc}^{(i-1)}_{n+1}  
	+ \jac{f}{y} \tilde {\bfg y}_{\rm int} ) \, .
\end{align}
%

We proceed as follows:  For $i=1$, substitution of \eqref{eq:adams:pred:lin} in \eqref{eq:adams:corr:lin} allows eliminating $\tilde {\bfg \xpc}^{(0)}_{n+1}$
and obtain $\tilde {\bfg \xpc}^{(1)}_{n+1}$.  
Subsequently, $\tilde {\bfg \xpc}^{(1)}_{n+1}$
is substituted in \eqref{eq:adams:corr:lin}, and so on $\forall i<r$, until we obtain $\tilde {\bfg \xpc}^{(r)}_{n+1}$.  Taking into account
\eqref{eq:corr:xt:lin}, we obtain the following form:
\begin{equation}
\begin{aligned}\label{eq:adams:xt:lin:2}
\tilde {\bfg x}_{n+1} &= 
\bfb R_1 \tilde {\bfg x}_{n} +
\bfb Q_1 \tilde {\bfg y}_{n} +
\bfb R_2 \tilde {\bfg x}_{n-1} +
\bfb Q_2 \tilde {\bfg y}_{n-1} 
\\
& + 
\ldots
+
\bfb R_{k} \tilde {\bfg x}_{n-k+1} +
\bfb Q_{k} \tilde {\bfg y}_{n-k+1} \, ,
\end{aligned}
\end{equation}
where $\bfb R_j$, $\bfb Q_j$ are proper 
coefficient matrices.  Considering \eqref{eq:pc:g:lin} and using the notation $\xs_n = (\tilde{\bfg x}_n, \tilde{\bfg y}_n)$, \eqref{eq:adams:xt:lin:2}, \eqref{eq:pc:g:lin} give:
\begin{equation}
\label{eq:adams:lin:aug}
\bfb K \xs_{n+1} = \bfb H_1 \xs_{n} + \bfb H_2 \xs_{n-1} + \ldots + \bfb H_{k} \xs_{n-k+1} \, ,
\end{equation}
\begin{equation}
\nonumber
\hspace{-8mm} 
{ \text{where}}
\hspace{4mm}
\bfb K =
\begin{bmatrix}
\bfb I   & \bfb 0 \\
\jac{g}{x} & \jac{g}{y} \\
\end{bmatrix}  , \ 
\bfb H_j =
\begin{bmatrix}
\bfb R_j & \bfb Q_j \\
\bfg 0   & \bfg 0 \\
\end{bmatrix}  , \
j \geq 1 \, .
\end{equation}
Then, by setting 
$\bfb y_{n+1}^{[0]} = \bfb x_{n+1}$,
$\bfb y_{n+1}^{[1]} = \bfb x_{n}$, 
$\ldots$ ,
$\bfb y_{n+1}^{[k-1]} = \bfb x_{n-k+2}$,
and
$\bfb y_{n}^{[0]} = \bfb x_{n} = \bfb y_{n+1}^{[1]}$,
$\bfb y_{n}^{[1]} = \bfb x_{n-1} = \bfb y_{n+1}^{[2]}$,
$\bfb y_{n}^{[k-1]} = \bfb x_{n-k+1}$.
Then, \eqref{eq:adams:lin:aug} can be equivalently written as:
\begin{equation}
\nonumber
\bfb H_{k} \bfb y_{n}^{[k-1]} = 
\bfb K \bfb y_{n+1}^{[0]}  
-\bfb H_1 \bfb y_{n+1}^{[1]} -
\bfb H_2 \bfb y_{n+1}^{[2]}
- \ldots -
\bfb H_{k-1} \bfb y_{n+1}^{[k-1]}
 \, .
\end{equation}
Collecting the relationships above in a matrix form, we arrive at \eqref{eq:difference}, with
$\bfb y_n = (\bfb y_n^{[0]}, \bfb y_n^{[1]}, \ldots, \bfb y_n^{[k-1]})$, and:
\begin{equation}
\nonumber
\bfg {\mathcal{E}} =
\begin{bmatrix}
\bfg 0   & \bfb I  \\
\bfb K & -\bfb H \\
\end{bmatrix}  , \ 
\bfg {\mathcal{A}} =
\begin{bmatrix}
\bfb I & \bfg 0 \\
\bfg 0 & \bfb H_{k} \\
\end{bmatrix}  ,
\
\bfb H = [\bfb H_1 \ \bfb H_2 \ \ldots \ \bfb H_{k-1}] \, .
\end{equation}

\bibliographystyle{IEEEtran}
\bibliography{references}

  \begin{biography}
    [{\includegraphics[width=1in,height=1.25in,clip,keepaspectratio]{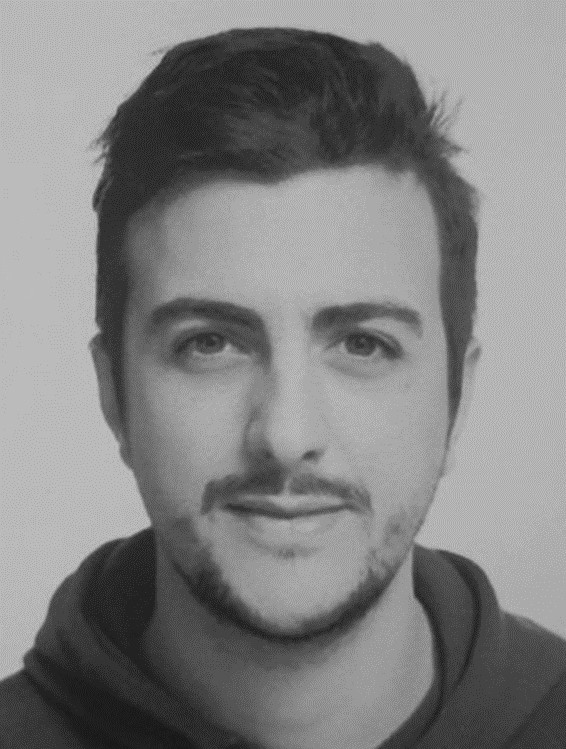}}]
    {Georgios Tzounas} (M'21) received the Diploma (M.E.) in Electrical and Computer Engineering from the National Technical Univ. of Athens, Greece, in 2017, and the Ph.D. 
    from Univ. College Dublin (UCD), Ireland, in 2021. In Jan.-Apr.~2020, he was a visiting researcher at Northeastern Univ., Boston, MA. 
     From Oct. 2020 to May 2022, he was a Senior Researcher at UCD.  He is currently a Postdoctoral Researcher at ETH
    Z{\"u}rich, working with the Power Systems Laboratory and the NCCR Automation project.    His research interests include modelling, stability analysis and control of power systems.
  \end{biography}
  
  \begin{biography}
    [{\includegraphics[width=1in,height=1.25in,clip,keepaspectratio]{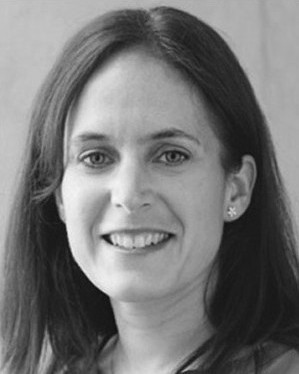}}]
    {Gabriela Hug} (SM'14) was born in Baden, Switzerland.  She received the M.Sc. degree in electrical engineering and the Ph.D. degree from the Swiss Federal Institute of Technology (ETH), Zürich, Switzerland, in 2004 and 2008, respectively.  After the Ph.D. degree, she worked with the Special Studies Group of Hydro One, Toronto, ON, Canada, and from 2009 to 2015, she was an Assistant Professor with Carnegie Mellon University, Pittsburgh, PA, USA.  She is currently a Professor with the Power Systems Laboratory, ETH Zurich.  Her research is dedicated to control and optimization of electric power systems.
  \end{biography}
  
\vfill

\end{document}